\magnification =1090 \hsize 16truecm \vsize 22.5truecm \hfuzz =0pt

\hfuzz =0pt
\scrollmode
\def\qed{\hbox{\vrule height 7pt depth 0pt width 7pt}}
\def\cqfd{\hfill\penalty 500\kern 10pt\qed\medbreak}

\def \l{{\lambda}}
\def \a{{\alpha}}
\def \b{{\beta}}
\def \s{{\sigma}}

\def \o{{\omega}}
\def \O{{\Omega}}
\def \d{{\delta}}

\def \R{{\bf R}}
\def \Z{{\bf Z}}
\def \N{{\bf N}}

\def \D{{\Delta}}
\def \E{{\bf E}\,}

\def \A{{\cal A}}

 \def \P{{\bf P}}

\def \qq{{\qquad}}
 \def \noi{{\noindent}}
 
 \def \p{{\varphi}}
\def \e{{\varepsilon}}
 \def \t{{\theta}}

\font\ph=cmcsc10  at  10 pt
\font\iit =cmmi7 at 7pt

\font\ph=cmcsc10  at  10 pt

 \font\gem= cmcsc10 at 10 pt
\font\gum= cmbx10 at 12 pt
\font\gem= cmbx10 at 10  pt

\rm
\def\ddate {\ifcase\month\or January\or
February\or March\or April\or May\or June\or July\or
August\or September\or October\or November\or December\fi\ {\the\day},
{\the\year}}
 
 \centerline{\gum Divisors of Bernoulli sums}
     \vskip 0,5 cm
  \centerline{ Michel {\ph Weber}}
\vskip 0,5 cm
  {\leftskip =2cm \rightskip=2cm   \noindent {\bf Abstract:\rm\   Let $ B_n= b_1+\ldots +b_n$, $n\ge 1$ where $ b_1,b_2,\ldots $ are 
independent Bernoulli random variables. In relation with the divisor problem, we evaluate the   almost sure asymptotic
order  of the   sums  $ \sum_{n=1}^N d_{\t, {\cal D}}(B_n )$, where $d_{\t, {\cal D}}(B_n )  = \#\{  
d\in {\cal D},     d\le   n^\t \! : \! d|B_n\}$  and   ${\cal D}$ is a  sequence of positive integers.
 }  
   \par }

\footnote {}{{\iit Date}\sevenrm : on \ddate\par \vskip -2pt   {\iit AMS} \ {\iit Subject}\ {\iit Classification} 2000:  Primary 11M06,
 Secondary 11M99, 11A25, 60G50. \par \vskip -2pt
  {\iit Keywords}:     Bernoulli random variables, iid, Theta functions,  
 divisors.\par } 
   \vskip 0,4cm
\medskip\noi\centerline{\gum 1. Introduction and main results}
\medskip\medskip\noi We begin with discussing a simple problem.  Let $d(n)=\#\{y :y|n\}$ be the divisor function, and consider also the
prime divisor function $\o(n)=\#\{p \, {\rm prime}:p|n\}$. Let $\b= \{\b_i , i\ge 1\}$ be a Bernoulli
sequence  and denote  
$B_n  =
\b_1+\ldots +\b_n$,
$n=1,2,\ldots$   the sequence of associated partial sums. Let   $(\O, \A,
\P)$ be  the underlying  basic probability
space. 
 It is natural to consider
the sequence of sums
 $$\sum_{n=1}^N d(B_n), \qq N=1,2,\ldots$$
and ask for which nondecreasing functions (optimal  if possible) $\Phi_1, \Phi_2:\N\to \R^+$, the following almost sure asymptotic
behavior can be established:   
$$ \sum_{n=1}^N d(B_n)\buildrel{a.s.}\over{=} \Phi_1(N) + {\cal O}(\Phi_2(N)). \eqno(1.1)$$
A similar question can naturally be raised relatively to the  prime divisor function $\o(n) $. Before giving an
arithmetical motivation for studying this problem, we would like to begin with a first necessary comment. A result of this sort cannot   be
obtained from the knewledge of the similar known result  for the deterministic sums:
$$\sum_{n=1}^N d( n)= \p_1(N) + {\cal O}(\p_2(N)).\eqno(1.2)$$
Indeed, the two sums    $\sum_{n=1}^N d(B_n)$  and $\sum_{n=1}^N d( n)$ are different, the first contains terms which appear with some
multiplicity, the order of this one can be bounded by $c\log n$, $n$ large, almost surely. And the natural idea to use  the law of the
iterated logarithm ($\forall\e>0$)
$$\big|B_n-{n\over 2}\big|\le \sqrt {(1+\e)n\log \log n}\quad n\ {\rm ultimately,\qq almost\ surely.} \eqno(1.3)$$ 
in order to exploit   (1.2) will give a less precise result than the one expected in (1.1). 
The law of the iterated
logarithm in (1.3) involves intervals of integers of the type $[m, m+C\sqrt {m\log \log m}]$. The study of the size of the divisor
function $d(n)$ for $n$ varying in intervals of this type or $[m, m+C\sqrt {m }]$, is usually known as the model of small intervals. 
Here the problem considered involves another (probabilistic) model, the one generated by the complete Bernoulli random walk; and so  
     relies upon the study of the asymptotic evolution of the system 
$$ \big\{ d(B_n), n\ge 1\big\}.$$
We are thus led to a probabilistic question  and the first natural object of investigation should consist of making a complete second
order theory of the above system, more precisely, a study of the correlation  
$$\E\, d(B_n) d(B_m)- \E\, d(B_n)\E\, d(B_m) \qq m>n, \eqno(1.4)$$
and preliminarily of (denoting $\chi$ the indicator function)
$$\E\, \chi(d|B_n) \chi(\d|B_m)- \P\{ d|B_n\}\P\{ \d|B_m\} \qq m>n. \eqno(1.5)$$
  Such a study turns up to depend, via the use of characteristic functions, on a careful analysis on the
circle of some naturally related $ cosine$ sums. 
 
  An arithmetical motivation to the study of   (1.1) can be easily provided.     Unlike to what happens
for the sum $\sum_{n=1}^N \o( n)$, where very accurate estimations of the order  of magnitude are known, the similar problem for the sum
$\sum_{n=1}^N d( n)$ contains a yet unsolved and certainly quite hard conjecture. For the comments  we shall now make, we refer to the
paper   [IM] of    Iwaniec and   Mozzochi. More
precisely, let
$$D(x) =\sum_{1\le n\le x}d(n) , \qq \D(x)= D(x)-x\log x -(2\gamma-1) x,\eqno(1.6)$$
where $\gamma$  is Euler's constant, and let   $\t_0$ be the smallest value of $\t$  such that
 $\D(x)\ll x^{\t+\e}$. 
It is conjectured since the papers of   Hardy (1916) and Ingham (1940)  that 
 $\t_0=1/4$. The best known result is (see the quoted paper)
 $ \t= 7/22\approx0,31181818\ldots $
 The study of the correlation problem described before should allow to obtain as corollary, via a suitable form of G\'al-Koksma's
Theorems, a result of the type 
 $$ \sum_{1\le n\le x}d(B_n)\buildrel{a.s.}\over{=}\sum_{1\le n\le x}\E \,d(B_n)+ {\cal O}_\e\Big(\big( \sum_{1\le n\le x}\E
\,d(B_n)\big)^{1/2+\e}\Big).
\eqno(1.7) $$
  Similarly the study of the quadruple correlation 
$$\E\, \big(d(B_n)-\E\, d(B_n)\big)\big(d(B_m)-\E\, d(B_m)\big)\big(d(B_p)-\E\, d(B_p)\big)\big(d(B_q)-\E\, d(B_q)\big)    \eqno(1.8)
$$  would provide by means of the same convergence criteria a result of the type:
 $$ \sum_{1\le n\le x}d(B_n)\buildrel{a.s.}\over{=}\sum_{1\le n\le x}\E \,d(B_n)+ {\cal O}_\e\Big(\big( \sum_{1\le n\le x}\E
\,d(B_n)\big)^{1/4+\e}\Big).
\eqno(1.9) $$
 \noi   That the above could be derived from (1.8) is already a remarkable fact; and it is clear from the very form of this result,
also in the light  of the
Dirichlet conjecture, that it would be  of considerable  interest  in succeeding to prove (1.9). Even a weaker form of it involving the
truncated divisor function
  $d_\t(m)=\#\{y \le m^\t\!:\!y|m\}   $, ( $ \theta\le  1/2$)
   would be also remarkable.  
 In this paper, we explore the correlation problem of order two (thus related to (1.4), (1.5)) and obtain almost sure results towards
(1.1).     
    Consider  as well   divisor functions defined with respect to a prescribed set  ${\cal D}$ of integers.  
  Define   
$$   d_{  {\cal D}}( n )  = d_{\t, {\cal D}}( n )  =\sum_{  
d\in {\cal D}  \atop   d\le   n^\t  } {\bf 1}_{\{d| n\}}.\eqno(1.10) $$
 The order of magnitude of this generalized divisor function (with $\t=1/2$) was examined in [We4] with the help of a
randomization argument (see related works in [BW], [We1-4]). 
 Let also $\eta>0$.   For divisors related to Bernoulli sums    we put similarly
$$  d_{ \eta, {\cal D}}(B_n )  =\sum_{ d<\eta \sqrt{n\over \log n}  \atop  d\in {\cal D} }\qq \qq  d_{\t, {\cal D}}(B_n )  =
 \sum_{ d\le  
n^\t  
   \atop    d\in {\cal D} } {\bf 1}_{\{d|B_n\}}. $$
 
The main purpose of the present study will consist of establishing limit
theorems for the sums 
$$\sum_{n=1}^N   d_{\t, {\cal D}}(B_n ) ,\qq \qq \sum_{  n\le N   \atop n\in{\cal N}}   d_{  \eta, {\cal D}} (B_{n} )\eqno(1.11)   $$
where ${\cal N} $ is an increasing sequence of positive integers. The difficult case is   ${\cal N}=\N$  and our results will be
then less precise than  in the subsequence case. When  ${\cal N}=\{\nu_k, k\ge 1\}$ satisfies for some 
$\rho>0$   the   growth
condition   
$$\nu_{k+1}-\nu_k\ge C\nu_k^\rho,\qq \qq k\ge 1. \eqno({\cal G}_\rho) $$
the second sum in (1.11) is controlable for $\eta\le \eta_\rho$, where $\eta_\rho$ depends on $\rho$ only, {\it no matter} ${\cal D}$ is.
When ${\cal N}=\N$,  restrictions on the range of $\t$ arise.  Naturally these estimates  rely upon  ${\cal D}$ and ${\cal N}$.  More
precisely  we show that these sums are almost surely  asymptotically comparable to their respective  (computable) means
$$\sum_{n=1}^N  \E d_{\t, {\cal D}}(B_n ) ,\qq \qq \sum_{  n\le N   \atop n\in{\cal N}}  \E d_{  \eta, {\cal D}} (B_{n} )   $$
and give an already sharp (although not optimal) estimate of the approximation rate. In the above case 
  by using (1.16) next (1.17) below, we get
$$ M_{\eta, {\cal N},  {\cal D}}(N):=\sum_{  n\le N   \atop n\in{\cal N}}  \E d_{  \eta, {\cal D}} (B_{n} )= \sum_{  n\le N\, ,\, d<\eta
\sqrt{n\over \log n}   \atop  n\in{\cal N}, d\in {\cal
D}         }{1\over d}  +{\cal O} (1).\eqno(1.12)$$
 And  
$$ M_{\t,  {\cal D}}(N):=\sum_{  n\le N   }  \E d_{  \t, {\cal D}} (B_{n} )=  \sum_{ n\le N\, ,\, d\le  
n^\t  
   \atop    d\in {\cal D} } \P{\{d|B_n\}}=\sum_{  n\le N\, ,\, d<n^\t\atop    d\in {\cal
D}         }{1\over d}  +{\cal O} (1) .\eqno(1.13)$$

Before stating the results  we shall first  comment more on correlation problem. This is a central question in the paper and section  2 
is entirely  devoted to its study.  The crucial point concerns the obtention of sharp 
 estimates for the   correlation function
$$ {\bf \Delta} \big( (d,n), (\d, m)\big)= \P\big\{ d|B_n\, ,\,   \d|B_m\big\}-\P \{ d|B_n \}\P \{    
\d|B_m \}, $$ 
and also for  the probability
$$   \P\big\{ d|B_n\, ,\,   \d|B_m\big\}. $$ 

\noi There are two cases of very different nature and unequal difficulty: the weakly dependent case ($n+n^c\le m$) which is relatively
easy to treat,  
 and
the dependent case ($n\le m\le n+n^c $). Here $c$ is some very small positive real. In the weakly dependent case, there is a constant
$C$ depending on $c$ such that for all $\eta$
  sufficiently small   and   $n$ large enough (see (3.3))
$$ \line{$  \displaystyle{\sup_{d<\eta  \sqrt{n\over \log n}  \atop \d<\eta \sqrt{m\over\log m} 
}}\big|{\bf \Delta}
\big( (d,n), (\d, m)\big)\big|
 \le 
    \left\{\matrix{    C  n  ({  \log (m-n)\over m-n})^{ 1/c} 
     \qq\ \   &{\rm if}\   n+n^c \le m\le 2n ,
\cr &\cr    C  
n({  \log n\over n})^{{1/ 2c}}({ \log (m-n)\over
m-n})^{1/2} &      {\rm if}\    m\ge 2n .\qq\qq\cr
 }  \right.  \hfill$}  $$

The dependent case  is the difficult
case and the only way we found, after having tried others,  to bound efficiently
${\bf
\Delta} \big( (d,n), (\d, m)\big)$ was, to start with  ${\bf
\Delta} \big( (d,n), (\d, m)\big)  \le    \P\big\{ d|B_n\, ,\,   \d|B_m\big\}$,   next to compare $\P\big\{ d|B_n\, ,\,   \d|B_m\big\}$ with $  \P\big\{
d\d|B_n  B_m\big\} $, and estimate  the probability
 $   \P\big\{ D|B_n  B_m\big\} $. This is, however, not a simple task and involves truly number theoritical arguments.  Exponentials of
second order arise (in (2.28))  for which we used   Sark\"osy's estimate      (Lemma 2.13). And the
multiplicative functions
$\rho_k(D)=\#\big\{1\le r\le D:   D| r^2+kr \big\}$, $k\le m-n$,  play  a central role when $n$ becomes large. As a consequence of a
sharper result (Proposition 2.10) , we show in section 2 that 
   $$\P \{ 
d|B_n\, ,\, \d|B_m \}\le  C   (m-n){ 2^{\o(d\d) } \over d\d
}+ C_\e  {     (d\d)  ^{(1 + \e)/2}\over  \sqrt n } .$$
Although we are convinced that this
bound is quite sharp, we are less sure that it fully describes  what happens for ${\bf
\Delta} \big( (d,n), (\d, m)\big)$ in the dependent range $(n, n+n^c)$, and must say that we have no alternative clue at the present
time. 
 \smallskip We can now  state our main results 
\medskip\noi{\gem Theorem 1.1.} {\it Let $0<\t<1/6$.    Then 
for any
$\e>0$}, $$  \sum_{n=1}^N   d_{\t, {\cal D}}(B_n )  \buildrel{a.s.}\over{=}M_{\t,  {\cal D}}(N) +{\cal O}_\e\big( M^{1/2+\e}_{ \t, {\cal
D}}(N)
\big) 
   .
$$ 
  \medskip\noi{\gem Theorem
1.2.} {\it Let ${\cal N} $ be   satisfying the growth condition $({\cal G}_\rho)$ for some 
$\rho>0$.   Put
  $M_N = \sum_{  n\le N   \atop n\in{\cal N}}
      \log^4 n.$
  Then there exists  $\eta_\rho$ depending on $\rho$ only, such that for $\eta\le \eta_\rho$    
$$  \sum_{  n\le N   \atop n\in{\cal N}}   d_{  \eta, {\cal D}} (B_{n} )  \buildrel{a.s.}\over{=}M_{\eta, {\cal N},  {\cal D}}(N) +{\cal
O}_\e\Big(  M_N^{1/2} 
\log^{3/2+\e}M_N\Big)    .
 $$ 
   And if ${\cal N}$ grows at most polynomially, then for some constant $b_0$,} 
 $$   \sum_{  n\le N   \atop n\in{\cal N}}   d_{  \eta, {\cal D}} (B_{n} )   \buildrel{a.s.}\over{=}M_{\eta, {\cal N},  {\cal D}}(N)
+{\cal O}_\e \Big( 
 M_{\eta, {\cal N},  {\cal D}}(N)^{1/2}
\big(\log  M_{\eta, {\cal N},  {\cal D}}(N)\big)^{b_0+\e}\Big)    .
 $$ 
 From the proof given in section 4  follows that $b_0 >7/2 $  suffices, but this value is certainly far from being optimal.
    Getting an
optimal rate of approximation appears as a certainly difficult and quite challenging question. It is also clear from the proofs of the
results, we shall give in the next sections,  that the error term is however improvable under    additional conditions on the sequence
${\cal D}$. Relevant conditions are for instance of the  type  
$$(a) \quad \#\big\{{\cal D}\cap [1,N]\big\}={\cal O}(N^\tau) \qq{\rm or}\qq (b) \quad \#\big\{{\cal D}\cap [1,N]\big\}={\cal
O}(\log^\eta N), \eqno( {\cal C} )$$ 
for some $0<\tau<1$ or  $\eta>0$. But this aspect of the
problem is not   considered in the present study. 

\smallskip For proving these results, our essential task will be   to bound efficiently the increments
 $\E\big(\sum_{i\le n\le j} H_n\big)^2$, 
where
 $  H_n=\sum_{d\le n^\t, d\in {\cal D}}\big({\bf 1}_{\{ d|B_n \}}-\P \{ d|B_n \}\big)   $,
  and  it is clear that it suffices   to bound
$
{\bf \Delta}   $ instead of its absolute value.  
   Some already existing results on the value distribution of Bernoulli sums will be incorporated   into the proofs. We briefly
recall them. Consider the elliptic Theta function 
$$\displaystyle{ {\bf \Theta} (d,m)  =  \sum_{\ell\in \Z} e^{im\pi{\ell\over   d }-{m\pi^2\ell^2\over 2 d^2}}. 
   }$$
In [We3]  (Theorem II),    the following uniform estimate is established: 
 $$\sup_{2\le d\le n}\Big|\P\big\{d|B_n\big\}- {\Theta(d,n)\over d}  \Big|= {\cal O}\big((\log n)^{5/2}n^{-3/2}\big).\eqno(1.14)  $$
 Here and throughout the whole paper,  $C$ will denote some absolute constant, which may change of value at each occurence. It is
easily seen  that 
$$ 
 \big|{\Theta (d,n)\over d}-{1\over d}\big| \le  \cases{   { C\over d}
 e^{ - {n \pi^2 \over 2d^2}}   & \qq if\ $d\le \sqrt n$,\cr&\cr
       {C \over\sqrt n}   & \qq if\ $ \sqrt n\le d\le n$.\cr
}   $$ 
 Therefore $$ 
  \big|\P\big\{d|B_n\big\}- {1\over d}  \big|\le  \cases{  C\Big((\log n)^{5/2}n^{-3/2}+ { 1\over d}
 e^{ - {n \pi^2 \over 2d^2}} \Big)   & \qq if $d\le \sqrt n$,\cr&\cr
      {C \over\sqrt n}   & \qq if $ \sqrt n\le d\le n$.\cr
}  \eqno(1.15)$$
 Further for any $\a>0$  
$$ \sup_{d<  \pi   \sqrt{   n \over 2\a\log n}}\big|\P\big\{  d|   B_{ n} \big\}-{1\over d} 
\big|= {\cal O}_\e\big(n^{-\a+\e }\big),\qq \quad (\forall \e>0). \eqno(1.16)$$  
 and for any $0<\rho<1 $, 
$$ \sup_{d<  (\pi/\sqrt 2) n^{(1-\rho)/2} }\big|\P\big\{  d|   B_{ n} \big\}-{1\over d} 
\big|= {\cal O}_\e\big(e^{-(1-\e) n^\rho}\big),\qq \quad (\forall 0<\e<1). \eqno(1.17)$$  
 Estimate  (1.17) exhibits a dramatic variation of the
uniform speed of convergence  of $\P\big\{d|B_n\big\} $ to its limit ${1/ d}$, when switching from the case $d\le    
n^{1/2} $ to the case $d\le    
n^{\t} $,
 $\t<1/2$.   
 It follows   that $\lim_{n\to \infty}\P\big\{d|B_n\big\}= {1/ d}  $, and
$$  \big|\P\big\{d|B_n\big\}- {1\over d}  \big| \le C   {  d\over n},
\qq \qquad {\rm if}\ 2\le d\le \sqrt n.\eqno(1.18)$$
 In particular $$ \sup_{2\le d\le \sqrt n}d\P\big\{d|B_n\big\} \le        C  .\eqno(1.19)$$


\font\exposubsectionfont=cmbx10 at 9 pt

  \bigskip\noi\centerline{\gum 2. Second order theory of  (I{\lower 2,5pt\hbox{(\exposubsectionfont d$  |$B{\lower
2,5pt\hbox{\exposubsectionfont n}})}} -P(d$  |$B{\lower 2,5pt\hbox{\exposubsectionfont n}})){\lower 2,5pt\hbox{\exposubsectionfont n}}}
 \medskip\noi   Our starting point
is   the   formula  
 $u\d_{u |B_{n }}=\sum_{j=0}^{u-1} e^{2i\pi  {j \over u}B_{n }}$,  
from which we deduce  after integration
$$\P \{ u|B_n \}={1\over u}\sum_{j=0}^{u-1}\E e^{2i\pi  {j \over u}B_{n }}=  {1\over u}\sum_{j=0}^{u-1} \Big({e^{2i\pi n{j\over u}}
+1 \over 2}\Big)^{n}=  {1\over u}\sum_{j=0}^{u-1} e^{i\pi n{j\over u}} \cos^{n} {
\pi j\over u}  .\eqno(2.1)$$
Thereby
$$\P \{ d|B_n \}\P \{ \d|B_m \} ={1\over d\d}\sum_{j=0}^{d-1}\sum_{h=0}^{\d-1} e^{i\pi  ( n{j\over d}+ m {h\over \d} )} \cos^{n} {
\pi j\over d}     \cos^{m} {
\pi h\over \d}. \eqno(2.2) $$
Let $m\ge n$. Similarly
$$\eqalign{ \P\big\{ d|B_n\, ,\,   \d|B_m\big\}    & ={1\over d\d}\E \Big\{ \sum_{j=0}^{d-1}\sum_{h=0}^{\d-1} e^{2i\pi  (  {j\over
d}B_n+   {h\over
\d}B_m )}\Big\} ={1\over d\d}  \sum_{j=0}^{d-1}\sum_{h=0}^{\d-1} \E e^{2i\pi  (  {j\over d} +   {h\over \d})B_n } \E e^{2i\pi   {h\over
\d} B_{m-n} }
\cr& = {1\over d\d}  \sum_{j=0}^{d-1}\sum_{h=0}^{\d-1}  e^{ i\pi  (  {j\over d} +   {h\over \d}) n }\cos^n{ \pi  (  {j\over d} +  
{h\over \d})   }  e^{ i\pi   {h\over
\d} (m-n) }\cos^{m-n}{ \pi   {h\over
\d}  }
\cr& = {1\over d\d}  \sum_{j=0}^{d-1}\sum_{h=0}^{\d-1}  e^{ i\pi  (  {j\over d}n +   {h\over \d}m) }\cos^n{ \pi  (  {j\over d} +  
{h\over \d})   }   \cos^{m-n}{ \pi   {h\over
\d} }\cr& = {1\over d\d}  \sum_{j=1}^{d-1}\sum_{h=1}^{\d-1}  e^{ i\pi  (  {j\over d}n +   {h\over \d}m) }\cos^n{ \pi  (  {j\over d} +  
{h\over \d})   }   \cos^{m-n}{ \pi   {h\over
\d} } \cr &\ \  +{1\over d\d}\Big(   \sum_{h=0}^{\d-1}  e^{ i\pi      {h\over \d}m  }   \cos^{m }{ \pi   {h\over
\d} }+   \sum_{j=1}^{d-1}   e^{ i\pi     {j\over d}n   }\cos^n{ \pi    {j\over d}   
   }  \Big)
\cr& = {1\over d\d}  \sum_{j=1}^{d-1}\sum_{h=1}^{\d-1}  e^{ i\pi  (  {j\over d}n +   {h\over \d}m) }\cos^n{ \pi  (  {j\over d} +  
{h\over \d})   }   \cos^{m-n}{ \pi   {h\over
\d} } \cr &\ \ +{\P \{ \d|B_m \}\over   d}+{\P \{ d|B_n \}\over  \d}   -{1\over d\d}  .\cr} $$
Therefore 
$$\displaylines{\P\big\{ d|B_n\, ,\,   \d|B_m\big\}\hfill\cr \hfill = {1\over d\d}  \sum_{j=1}^{d-1}\sum_{h=1}^{\d-1}  e^{ i\pi  ( 
{j\over d}n +   {h\over
\d}m) }\cos^n{ \pi  (  {j\over d} +   {h\over \d})   }   \cos^{m-n}{ \pi   {h\over
\d} }  +{\P \{ \d|B_m \}\over   d}+{\P \{ d|B_n \}\over  \d}   -{1\over d\d}. \qq (2.3) \cr}  $$
And
$$\eqalign{ {\bf \Delta} \big( (d,n),  (\d, m)\big)   & ={1\over d\d}  \sum_{j=0}^{d-1}\sum_{h=0}^{\d-1}  e^{ i\pi  (  {j\over d}n +  
{h\over
\d}m) }  \cos^{m-n}{ \pi   {h\over
\d} }\Big\{\cos^n{ \pi  (  {j\over d} +  
{h\over \d})   }-\cos^{n} {
\pi j\over d}     \cos^{n} {
\pi h\over \d}\Big\}\cr & 
={1\over d\d}  \sum_{j=1}^{d-1}\sum_{h=1}^{\d-1}  e^{ i\pi  (  {j\over d}n +   {h\over \d}m) }  \cos^{m-n}{ \pi   {h\over
\d} }\Big\{\cos^n{ \pi  (  {j\over d} +  
{h\over \d})   }-\cos^{n} {
\pi j\over d}     \cos^{n} {
\pi h\over \d}\Big\}. \cr}\eqno(2.4) $$ 
Here the summands with $j=0$ or $h=0$ do not contribute.

\medskip\noi\centerline{\bf 2.1. Reductions   via symmetries}
\medskip\noi We begin with the probability 
$$\P \{ u|B_n \} =  {1\over u}\sum_{\ell=0}^{u-1} e^{i\pi n{\ell\over u}} \cos^{n} {
\pi \ell\over u}  . $$
\smallskip\noi --- If $u$ is odd, say $u=2r+1$,  and $r+1\le \ell\le
2r $, write
$\ell= 2r-\l$. Then
$ 0\le \l\le r-1$ and 
$$ {\ell\over u}= {2r-\l\over 2r+1}= 1 -{  \l+1\over 2r+1}:=1 -{  l\over u} \qq{\rm and} \qq 1\le l\le r.$$
Further 
$$\displaylines{e^{i\pi n{\ell\over u}} \cos^{n} {
\pi \ell\over u} = e^{i\pi n   -i\pi n {  l\over u} } \cos^{n}  (\pi -\pi{  l\over u})= e^{  -i\pi n {  l\over u} } \cos^{n}   
 \pi{  l\over u}   .\cr}$$
 Thus 
$$\P \{ u|B_n \} =  {1\over u}\sum_{ |\ell|<u/2} e^{i\pi n{\ell\over u}} \cos^{n} {\pi \ell\over u}  =  {1\over u}+ {2\over u}\sum_{ 
1\le \ell <u/2} \cos(\pi n{\ell\over u}) \cos^{n} {
\pi \ell\over u}  .\eqno(2.5)$$
    
\noi --- If $u$ is even: $u=2r$, then $\ell$ varies between $1$ and $r-1$, next between $r+1$ and $2r-1$ with a median value $\ell= r$.
For this indice, we have $\cos^{ n}{    {\pi \ell\over u}   } = \cos^{ n}{    {\pi r\over
2r}   }=\cos^{ n}{    {\pi  \over
2 }   }=0$, and there is no contribution. If $r+1\le \ell\le 2r-1$, write $\ell=2r-1-b$. Then $0\le
b\le r-2$ and
$$  {  \ell\over
u}=  {  2r-1-b\over
2r} = 1-{  b+1\over
2r}:=1 -{  l\over  u} \qq{\rm with} \qq 1\le l\le r-1.$$  
Thus 
 $ e^{i\pi n{\ell\over u}} \cos^{n} {
\pi \ell\over u} = e^{i\pi n   -i\pi n{  l\over u} } \cos^{n}  (\pi -\pi{  l\over u})= e^{  -i\pi n {  l\over u} } \cos^{n}   
 \pi{  l\over u}     $, 
and here again we have (2.5).
 \bigskip
Now, we pass to  the probability $\P\big\{ d|B_n\, ,\,   \d|B_m\big\}$. Here also  we   operate  reductions
allowing to work in the first quadrant only. This is quite similar to the above. By (2.3), $\P\big\{ d|B_n\, ,\,   \d|B_m\big\}
={\bf \Psi}+{\bf \Phi}$, where
$$\eqalign{{\bf \Psi} ={\bf \Psi}\big( (d,n), (\d, m)\big)&= {1\over d\d}  \sum_{j=1}^{d-1}\sum_{h=1}^{\d-1}  e^{ i\pi  (  {j\over d}n +  
{h\over
\d}m) }\cos^n{ \pi  (  {j\over d} +   {h\over \d})   }   \cos^{m-n}{ \pi   {h\over
\d} }  , \cr {\bf \Phi} ={\bf \Phi}\big( (d,n), (\d, m)\big)&= {\P \{ \d|B_m \}\over   d}+{\P \{ d|B_n \}\over  \d}   -{1\over d\d}.
\cr}\eqno(2.6) $$
  We shall thus   be mainly concerned with the sum ${\bf \Psi}$.
 
\smallskip\noi --- If $\d$ is odd, say $\d=2q+1$,  and $q+1\le h\le
2q $, write
$h= 2q-b$. Then
$ 0\le b\le q-1$ and 
$$ {h\over \d}= {2q-b\over 2q+1}= 1 -{  b+1\over 2q+1}:=1 -{  \b\over \d} \qq{\rm and} \qq 1\le \b\le q.$$
The corresponding summand writes 
$$\displaylines{e^{ i\pi m} e^{ i\pi  (  {j\over d}n -    {\b\over \d}m) }  \cos^{m-n}{ (\pi-     {\pi\b\over
\d} )}  \cos^n{ (\pi+ \pi  (  {j\over d} -  
{\b\over \d}))   }  \hfill\cr\hfill  = e^{ i\pi  (  {j\over d}n -    {\b\over \d}m) } \cos^{m-n}{ (     {\pi\b\over
\d} )}  \cos^n{    \pi  (  {j\over d} -  
{\b\over \d})    } .\cr}$$
 Thus 
$$\displaylines{{\bf \Psi}  ={1\over d\d}  \sum_{j=1}^{d-1}\Bigg\{\sum_{h=1}^{q}  e^{ i\pi  (  {j\over d}n +   {h\over \d}m) }  \cos^{m-n}{
\pi   {h\over
\d} } \cos^n{ \pi  (  {j\over d} +  
{h\over \d})   } +   \hfill\cr\hfill \sum_{h=-q}^{-1}e^{ i\pi  (  {j\over d}n +    {h\over \d}m) } \cos^{m-n}{    {\pi h\over
\d}   }   \cos^n{    \pi  (  {j\over d} +  
{h\over \d})    } \Bigg\}
\cr\hfill ={1\over d\d}  \sum_{j=1}^{d-1} \sum_{1\le |h|<\d/2}  e^{ i\pi  (  {j\over d}n +   {h\over \d}m) } 
\cos^{m-n}{
\pi   {h\over
\d} } \cos^n{ \pi  (  {j\over d} +  
{h\over \d})   } . \qq \qq \qq \qq \qq \qq    \       \quad  \cr}$$

\noi --- If $\d$ is even: $\d=2q$, then $h$ varies between $1$ and $q-1$, next between $q+1$ and $2q-1$ with a median value $h= q$. In the
latter case, we have $\cos^{m-n}{    {\pi h\over
\d}   } = \cos^{m-n}{    {\pi q\over
2q}   }=\cos^{m-n}{    {\pi  \over
2 }   }=0$, and there is no contribution of this indice. If $q+1\le h\le 2q-1$, write $h=2q-1-b$. Then $0\le
b\le q-2$ and
$$  {  h\over
\d}=  {  2q-1-b\over
2q} = 1-{  b+1\over
2q}:=1 -{  \b\over  \d} \qq{\rm and} \qq 1\le \b\le q-1.$$  
The corresponding summand writes exactly as before:
$$e^{ i\pi m} e^{ i\pi  (  {j\over d}n -    {\b\over \d}m) }  \cos^{m-n}{ (\pi-     {\pi\b\over
\d} )} \cos^n{ (\pi+ \pi  (  {j\over d} -  
{\b\over \d}))   } $$
and we have 
$$\displaylines{{\bf \Psi}  ={1\over d\d}  \sum_{j=1}^{d-1}\Bigg\{\sum_{h=1}^{q-1}  e^{ i\pi  (  {j\over d}n +   {h\over \d}m) } 
\cos^{m-n}{
\pi   {h\over
\d} } \cos^n{ \pi  (  {j\over d} +  
{h\over \d})   }   \hfill\cr\hfill+  \sum_{h=-(q-1)}^{-1}e^{ i\pi  (  {j\over d}n +    {h\over \d}m) } \cos^{m-n}{    {\pi h\over
\d}   }  \cos^n{    \pi  (  {j\over d} +  
{h\over \d})    } \Bigg\}\cr\hfill ={1\over d\d}  \sum_{j=1}^{d-1} \sum_{1\le |h|<\d/2}  e^{ i\pi  (  {j\over d}n +   {h\over \d}m) } 
\cos^{m-n}{
\pi   {h\over
\d} } \cos^n{ \pi  (  {j\over d} +  
{h\over \d})   } . \qq \qq \qq \qq \qq   \ \quad     \qq  \cr}$$
 A similar remark can be made concerning the  subsum $ \sum_{j=1}^{d-1}$. We   display this point again to make the proof transparent. 
\smallskip --- If $d$ is odd: $d=2p+1$, and $j$ is between $p+1$ and $2p$, write $j=2p-b$, $0\le b\le p-1$. Then 
$${j\over d} ={2p-b\over 2p+1}=1-{b+1\over 2p+1}:= 1-{\b\over d}\qq{\rm and} \qq 1\le \b\le p.$$  
The corresponding summand writes
 $$ e^{i\pi n}e^{i\pi (-{\b n\over d}+ {hm\over \d})}\cos^{m-n} {\pi h\over \d}  \cos^n(\pi+\pi(-{\b\over d} + {h\over
\d}))    =  e^{i\pi (-{\b n\over d}+ {hm\over \d})}\cos^{m-n} {\pi h\over \d}  \cos^n \pi(-{\b\over d} + {h\over
\d})  . $$
And 
$$\displaylines{{\bf \Psi}  ={1\over d\d}  \sum_{1\le |h|<\d/2} \Bigg\{  \sum_{j=1}^{p} e^{ i\pi  (  {j\over d}n +   {h\over \d}m) } 
\cos^{m-n}{
\pi   {h\over
\d} } \cos^n{ \pi  (  {j\over d} +  
{h\over \d})   } \hfill\cr\hfill + \sum_{j=1}^{p}e^{i\pi (-{j n\over d}+ {hm\over \d})}\cos^{m-n} {\pi h\over \d}  \cos^n \pi(-{j\over d}
+ {h\over
\d}) \Bigg\}\cr \hfill
={1\over d\d}  \sum_{1\le |h|<\d/2}   \sum_{1\le |j|<d/2}  e^{ i\pi  (  {j\over d}n +   {h\over \d}m) }  \cos^{m-n}{ \pi  
{h\over
\d} } \cos^n{ \pi  (  {j\over d} +  
{h\over \d})   } .\qq\qq\qq\qq    \ \ \ (2.7)\cr}$$
 
\smallskip --- If $d$ is even: $d=2p $, we distinguish between $1\le j\le p-1$, $p+1\le j\le 2p-1$ and the median value $j=p$, which this
time contributes to the sum. When $p+1\le j\le 2p-1$,  write
$j=2p-b-1$,
$0\le b\le p-2$. Then 
$${j\over d} ={2p-b-1\over 2p }=1-{b+1\over 2p }:= 1-{\b\over d}\qq{\rm and} \qq 1\le \b\le p-1.$$  
Thus the corresponding summand writes exactly as in the previous case
$$ e^{i\pi (-{\b n\over d}+ {hm\over \d})}\cos^{m-n} {\pi h\over \d}  \cos^n \pi(-{\b\over d} + {h\over
\d})  ,$$
and 
$$\displaylines{{\bf \Psi}  ={1\over d\d}  \sum_{1\le |h|<\d/2} \Bigg\{  \sum_{j=1}^{p-1} e^{ i\pi  (  {j\over d}n +   {h\over \d}m) } 
\cos^{m-n}{
\pi   {h\over
\d} } \cos^n{ \pi  (  {j\over d} +  
{h\over \d})   } \hfill\cr\hfill+ \sum_{j=1}^{p-1}e^{i\pi (-{j n\over d}+ {hm\over \d})}\cos^{m-n} {\pi h\over \d}  \cos^n
\pi(-{j\over d} + {h\over
\d}) \Bigg\}\cr \hfill
={1\over d\d}  \sum_{1\le |h|<\d/2}   \sum_{1\le |j|<d/2}  e^{ i\pi  (  {j\over d}n +   {h\over \d}m) }  \cos^{m-n}{ \pi  
{h\over
\d} } \cos^n{ \pi  (  {j\over d} +  
{h\over \d})   }  \qq\qq\qq\qq \ \ \ \ \qquad   
\cr \hfill+{1\over d\d}   \sum_{1\le |h|<\d/2} e^{ i\pi  (  {n\over 2}  +   {hm\over \d} ) }  \cos^{m-n}{     {\pi h\over
\d} } .\cos^n{    (  {\pi \over 2} +  
{\pi h\over \d})   }. \cr}$$ 
We  therefore get  
 $$ {\bf \Psi} = {1\over d\d}  \sum_{1\le |h|<\d/2}   \sum_{1\le |j|<d/2}  e^{ i\pi  (  {j\over d}n +   {h\over \d}m) }  \cos^{m-n}{
\pi   {h\over
\d} } \cos^n{ \pi  (  {j\over d} +  
{h\over \d})   }+{\bf r},\eqno(2.8)$$
where $$\line{$ \qq\qq {\bf r} =
     \left\{\matrix{   0 
     \qq\ \   &\hbox{\rm if $d$ is odd} 
\cr &\cr  \displaystyle{{1\over d\d}   \sum_{1\le |h|<\d/2}} e^{ i\pi  (  {n\over 2}  +   {hm\over \d} ) }  \cos^{m-n}{     {\pi h\over
\d} } .\cos^n{    (  {\pi \over 2} +  
{\pi h\over \d})   }  &     \ \ \hbox{\rm if $d$ is even.} \cr
 }  \right.  \hfill$}\eqno(2.8a)$$
  \medskip\noi  Finally  as $\P\big\{ d|B_n\, ,\,   \d|B_m\big\}
= {\bf \Psi} +{\bf \Phi}$, we   obtained

\smallskip
{\bf First reduced form:}
$$\eqalign{\P\big\{ d|B_n\, ,\,   \d|B_m\big\}
&={1\over d\d}  \sum_{1\le |h|<\d/2}   \sum_{1\le |j|<d/2}  \cos    (  {\pi j\over d}n +   {\pi h\over \d}m)    \cos^{m-n}{ \pi  
{h\over
\d} } \cos^n{ \pi  (  {j\over d} 
  +  
{h\over \d})   }\cr &\quad +{\bf r}+{\P \{ \d|B_m \}\over   d}+{\P \{ d|B_n \}\over  \d}   -{1\over d\d}. \cr}\eqno(2.9) $$ 
Consequently  we   have to estimate     four sums of   type  
 $${\bf \Psi}_{\e, \eta}={1\over d\d} \sum_{ 1\le j< d/2\atop 1\le h<\d/2}  e^{ i\pi  (  \eta{j\over d}n + \e  {h\over \d}m) } 
\cos^{m-n}{ \pi   {h\over
\d} } \cos^n{ \pi  ( \eta {j\over d} + \e 
{h\over \d})   }   ,\eqno(2.10)$$
where $\e,\eta\in \{-1,+1\}$. 
Turning to ${\bf \Delta}$  we observe that 
$$ {\bf \Delta}  ={\bf \Psi}-{1\over d\d}  \sum_{j=1}^{d-1}\sum_{h=1}^{\d-1}  e^{ i\pi  (  {j\over d}n +  
{h\over
\d}m) } 
\cos^{m }{
\pi   {h\over
\d} }  \cos^{n} {
\pi j\over d}  = {\bf \Psi}-\big(\P \{ d|B_n \} -  {1\over d}\big)\big(\P \{ \d|B_m \} -  {1\over \d}\big).   \eqno(2.11)   $$
Since 
$$\eqalign{  \big(\P \{ d|B_n \} -  {1\over d}\big)\big(\P \{ \d|B_m \} -  {1\over \d}\big)&={1\over d\d} 
\sum_{j=1}^{d-1}\sum_{h=1}^{\d-1}  e^{ i\pi  (  {j\over d}n +   {h\over \d}m) }  \cos^{m }{ \pi   {h\over
\d} }  \cos^{n} {
\pi j\over d}\cr &
={1\over d\d}\sum_{1\le |h|<\d/2\atop1\le |j|<d/2}  e^{i\pi( n{j\over d}+m{h\over \d})} \cos^{n} {
\pi j\over d}  \cos^{m} {
\pi h\over \d}, \cr}  $$
we get

\smallskip
{\bf Second reduced form:}
 
$$ {\bf \Delta} \big( (d,n), (\d, m)\big)    
={1\over d\d}\sum_{1\le |h|<\d/2\atop1\le |j|<d/2}    e^{ i\pi  (  {j\over d}n +   {h\over \d}m) }\cos^{m-n}{ \pi  
{h\over
\d} }  \Big\{  \cos^n{ \pi  (  {j\over d} +  
{h\over \d})   }  -    \cos^{n} {
\pi j\over d}  \cos^{n} {
\pi h\over \d}\Big\} +{\bf r}. \eqno (2.12)   $$

\bigskip\noi\centerline{\bf 2.2. The weakly dependent case (${\bf m\ge n+ n^{c} } $)}

\bigskip\noi  In what follows   $ c\in ]0,1[$ is some fixed     small real.   We are indeed interested in results
valid for
$c$ arbitrary small.   By     (2.11)  we know that
$${\bf \Delta}= {\bf \Psi}\ -\big(\P \{ d|B_n \} -  {1\over d}\big)\big(\P \{ \d|B_m \} -  {1\over \d}\big),$$
 where the sum  ${\bf \Psi}={\bf \Psi} \big( (d,n), (\d, m)\big)$ is defined in (2.6). Fix also some reals $\a,\a'  $,  depending on   $c$, such that
$\a >\a' >\max(3/2, 1/c)$. We shall consider two subcases.

\bigskip\centerline{\bf Case:  ${\bf  n+ n^{c} \le m\le 2n.} $}
\medskip\noi

 We shall first
establish the following 
\medskip\noi {\gem Proposition 2.1.} {\it  There exist    constants $C $ and $n_0$  depending on  $c$, such that for $n\ge n_0$
and 
$   n+ n^{c} \le m\le 2n$  
$$\displaylines{ \big|{\bf \Psi} \big( (d,n), (\d, m)\big)\big| \le  C \Bigg\{   n({  \log (m-n)\over m-n})^{ 1/c}+
{1\over d\d}  \sum_{0< 
{\pi j\over  d}
\le 2 ( {2\a\log (m-n) \over m-n} )^{1/2}\atop 0< {\pi h\over \d}\le  ( {2\a\log
(m-n) \over m-n} )^{1/2}}     e^{-   {3n}  (  
{j\over d} -
 {h\over \d})^2    } 
 \Bigg\}. \cr}$$
In particular, if 
 $  
\max(d, \d)<  {\pi \over2\sqrt{2\a}   } ( {m-n\over \log (m-n)   } )^{1/2}$, 
   we have the simpler bound} 
$$\big|{\bf \Psi} \big( (d,n), (\d, m)\big)\big| \le  C     n({  \log (m-n)\over m-n})^{ 1/c}.$$
\medskip  The second assertion simply follows from the fact that the sum in the righthand side disappears. The proof of the above
proposition will result from   three lemmas. By (2.10), the sum
${\bf
\Psi}={\bf
\Psi}\big( (d,n), (\d, m)\big)$ is a sum of four subsums, which are all of the type 
 
$${\bf \Psi}_{\e, \eta}={1\over d\d} \sum_{ 1\le j< d/2\atop 1\le h<\d/2}  e^{ i\pi  (  \eta{j\over d}n + \e  {h\over \d}m) } 
\cos^{m-n}{ \pi   {h\over
\d} } \cos^n{ \pi  ( \eta {j\over d} + \e 
{h\over \d})   }   , $$
where $\e,\eta\in \{-1,+1\}$. We shall therefore estimate the sums ${\bf \Psi}_{\e, \eta}$ which, in what follows we will simply write
 ${\bf \Psi} $, when   no confusion. 

\medskip   Put for any integer
$v\ge 1$
 $$\p_v=   \big( {2\a\log v \over v}\big)^{1/2},
\qquad\  \tau_v= {\sin\p_v/2\over
\p_v /2} .\eqno(2.13)$$ 
  Let $n_0 $   be sufficiently large  so that for   $m-n\ge n_0$, 
$$  \tau_{m-n}\ge (\a^\prime/\a)^{1/2} .\eqno(2.14)$$    
 We distinguish between three cases: 
\smallskip \noi \qq\qq\qq\qq  --- (I) $\big( {2\a\log (m-n) \over m-n}\big)^{1/2}\le {\pi h\over \d} \le \pi/2 $.
\smallskip \noi \qq\qq\qq\qq  --- (II) $0< {\pi h\over \d}\le \big( {2\a\log (m-n) \over m-n}\big)^{1/2}$ and $0< 
{\pi j\over  d}
\le 2\big( {2\a\log (m-n) \over m-n}\big)^{1/2}$.
\smallskip \noi \qq\qq\qq\qq  --- (III) $0< {\pi h\over \d}\le \big( {2\a\log (m-n) \over m-n}\big)^{1/2}$ and $2\big( {2\a\log
(m-n) \over m-n}\big)^{1/2}\le {\pi j\over  d} \le \pi/2 $.

\medskip    {\bf Case I}. Put 
$${\bf \Psi}_{\bf 1}:= {1\over d\d}  \sum_{1\le j<{d\over 2}}\sum_{ ( {2\a\log (m-n) \over m-n} )^{1/2}\le {\pi h\over \d} \le
\pi/2}   e^{ i\pi  ( 
 {\eta j\over d}n +   {\e h\over \d}m) }  \cos^{m-n}{
\pi   {h\over
\d} } \cos^n{ \pi  (  {\eta j\over d} +  
{\e h\over \d})   }.$$
\noi{\gem Lemma 2.2.} {\it  There exist    constants $C $ and $n_0$ depending on  $c$  such that for $m,n$ such that $m-n\ge
n_0$,}
$$ |{\bf \Psi}_{\bf 1}| \le (m-n)^{- \a'}. $$  
  {\it Proof.} As  $|\cos {\pi h\over \d}|\le \cos \p_{m-n}$, we get since $\log u\le u-1$, $u>0$: for $m-n$ large enough, say $m-n\ge
n_0$ (so that $2\sin^2(\p_{m-n}/2)<1$)
 $$|\cos {\pi h\over \d}|^{m-n} \le  \cos^{m-n} \p_{m-n} =   e^{ (m-n)\log (1-2\sin^2(\p_{m-n}/2)})
\le e^{-2(m-n)\sin^2(\p_{m-n}/2)}  . $$
 But    
$$ 2(m-n) \sin^2(\p_{m-n}/2) = 2(m-n) (\p_{m-n}/2)^2\tau_{m-n}^2\ge     {\a'(m-n) \log m-n\over m-n}=  \a'  \log m-n 
    .$$ 
 Hence,  
$$\eqalign{|{\bf \Psi}_{\bf 1}|&\le {1\over d\d}  \sum_{1\le j<{d\over 2}}\sum_{ ( {2\a\log (m-n) \over m-n} )^{1/2}\le {\pi
h\over \d} \le \pi/2}     |
\cos {
\pi   {h\over
\d} }|^{m-n} \cr& \le {1\over  \d}  \sum_{ ( {2\a\log (m-n) \over m-n} )^{1/2}\le {\pi h\over \d} \le \pi/2}  \cos^{m-n} \p_{m-n}
\cr&\le \cos^{m-n}
\p_{m-n} 
 \le (m-n)^{- \a'}.\cr} $$  
 \cqfd

\medskip    {\bf Case II}.  Let 
$$\displaylines{{\bf \Psi}_{\bf 2} \hfill\cr\hfill  ={1\over d\d} \sum_{0<
{\pi j\over  d}
\le 2 ( {2\a\log (m-n) \over m-n} )^{1/2}}\sum_{0< {\pi h\over \d}\le  ( {2\a\log
(m-n) \over m-n} )^{1/2}}   e^{ i\pi  (  \eta{j\over d}n + \e  {h\over \d}m) } 
\cos^{m-n}{ \pi   {h\over
\d} }\ \cos^n{ \pi  ( \eta {j\over d} + \e 
{h\over \d})   }     .\cr}  $$

 We shall prove 
\medskip\noi {\gem Lemma 2.3.} {\it There exist    constants $C $ and $n_0$ depending on $c$ only such that for $\min(m-n, n)\ge n_0$  
} 
$$  |{\bf \Psi}_{\bf 2}
|\le  C  ({  \log (m-n)\over m-n})^{  1/c}+{1\over d\d}
  \sum_{0< 
{\pi j\over  d}
\le 2 ( {2\a\log (m-n) \over m-n} )^{1/2}\atop 0< {\pi h\over \d}\le  ( {2\a\log
(m-n) \over m-n} )^{1/2}}     e^{-   {3n}  (  
{j\over d} -
 {h\over \d})^2    }    . $$
The estimate only makes sense if $m> n+n^c$.

\medskip\noi {\it Proof.}   Using Mac-Laurin formula, for any positive integer $p$ there exists a polynomial $Q_p(x)=\sum_{s=1}^{p-1} a_sx^{2s}$, with $a_1=-1/2$,
$a_2=-1/12$,... and   constants $x_p, C_p$ depending on $p$ only, such that for $  |x|\le x_p$,
$$\big|\log \cos x -Q_p(x)\big|\le C_p |x|^{2p},$$
and 
$$\big|     Q_p(x)+x^2/2\big|\le C_p |x|^{4}.$$
Hence, for   $|x|\le x_p$, 
 $$      Q_p(x) \le   - x^{2}/3,$$
provided that $x_p$ is sufficiently small, which we do assume from now on. We select a integer $p$ so that 
$$ c> {1\over  p +1}. \eqno(2.15) $$ 
  Put $E(x)=e^{Q_p(x)}$. We shall compare the subsum 
 ${\bf \Psi}_{\bf 2} $
with
 
$$\displaylines{ {\bf \Psi}_{\bf 2}'\hfill\cr\hfill  ={1\over d\d} \sum_{0<
{\pi j\over  d}
\le 2 ( {2\a\log (m-n) \over m-n} )^{1/2}}\sum_{0< {\pi h\over \d}\le  ( {2\a\log
(m-n) \over m-n} )^{1/2}}    e^{ i\pi  (  \eta{j\over d}n + \e 
{h\over \d}m) } 
\cos^{m-n}{ \pi   {h\over
\d} } E^n(\pi  (  {\eta j\over d} +  
{\e h\over \d}))    . \cr} $$
 By using   the elementary inequality:  
$|e^u-e^v|\le |u-v| $ for $u,v\le 0$, we   observe that 
$$\big|\cos^n x -e^{ nQ_p(x)}\big|\le n\big|\log\cos  x - Q_p(x)\big|\le C_pn|x|^{2p}, $$
for   $|x|\le x_p$. We have $\pi | {\eta j\over d} +
 {\e h\over \d}|\le x_p$ once  $m-n$ is large enough.    Thus
$$  \Big|\cos^n{ \pi  (  {\eta j\over d} +  
{\e h\over \d})   }- E^n(\pi  (  {\eta j\over d} +  
{\e h\over \d})) \Big| \le  C_pn\big[ ({h\over \d})^{2p} +    ({j\over
d})^{2p} \big].  $$
  And so
$$\bigg|{1\over d\d} \sum_{0<
{\pi j\over  d}
\le 2 ( {2\a\log (m-n) \over m-n} )^{1/2}\atop 0< {\pi h\over \d}\le  ( {2\a\log
(m-n) \over m-n} )^{1/2}}    e^{ i\pi  (  \eta{j\over d}n +   \e{h\over \d}m) }  \cos^{m-n}{
\e   {\pi h\over
\d} }\Big\{\cos^n{ \pi  ( \eta {j\over d} +  
\e{h\over \d})   }-E^n(\pi  (  {\eta j\over d} +  
{\e h\over \d})) \Big\}\bigg|$$
$$\eqalign{& \le C {n\over d\d}  \sum_{0< 
{\pi j\over  d}
\le 2 ( {2\a\log (m-n) \over m-n} )^{1/2}\atop 0< {\pi h\over \d}\le  ( {2\a\log
(m-n) \over m-n} )^{1/2}}  
      \big[ ({h\over \d})^{2p} +    ({j\over
d})^{2p} \big] .\cr} $$ 

 Therefore
$$ \big|{\bf \Psi}_{\bf 2}-{\bf \Psi}'_{\bf 2}\big|\le  C   n\big({  \log (m-n)\over
m-n }\big)^{  p+1  } .\eqno(2.16)$$
 
\noi   
Further  since  
 $E(x)  \le e^{-    x ^2/3 }  $ for $|x|\le x_p$, we may also simply bound $|{\bf \Psi}_{\bf 2}'|$ for $n$ large enough as
follows  
$$\eqalign{  |{\bf \Psi}_{\bf 2}'|\le {1\over d\d} \sum_{0< 
{\pi j\over  d}
\le 2 ( {2\a\log (m-n) \over m-n} )^{1/2}\atop 0<{\pi h\over \d}\le  ( {2\a\log
(m-n) \over m-n} )^{1/2}}     e^{-  {3n} (  {\eta j\over d} +
 {\e h\over \d})^2 }  .
 \cr}\eqno(2.17)$$  
   It is clear that   these sums are bounded by   
$$   {1\over d\d} \sum_{0<
{\pi j\over  d}
\le 2 ( {2\a\log (m-n) \over m-n} )^{1/2}\atop 0<{\pi h\over \d}\le  ( {2\a\log
(m-n) \over m-n} )^{1/2}}  e^{-   {3n}  (  
{j\over d} -
 {h\over \d})^2    }     \eqno(2.18)$$
    With  estimates (2.16),  (2.17) and (2.18), we therefore get 
$$  |{\bf \Psi}_{\bf 2}
|\le  C \bigg\{n ({  \log (m-n)\over m-n})^{  p+1}+{1\over d\d}
  \sum_{0<
{\pi j\over  d}
\le 2 ( {2\a\log (m-n) \over m-n} )^{1/2}\atop 0< {\pi h\over \d}\le  ( {2\a\log
(m-n) \over m-n} )^{1/2}}     e^{-   {3n}  (  
{j\over d} -
 {h\over \d})^2    }\bigg\}  . $$
   
\noi    With $p$ chosen accordingly with (2.15), we have ${  \log (m-n)\over m-n})^{  p+1}\le {  \log (m-n)\over
m-n})^{ 1/c}$; so that if $n$ is sufficiently large, say $n\ge n_0$ where $n_0$ depends on $c$, $\a$ and $m\ge n+n^c $, we get
  $$  |{\bf \Psi}_{\bf 2}
|\le  C \bigg\{ n({  \log (m-n)\over m-n})^{   1/c}+{1\over d\d}
  \sum_{0< 
{\pi j\over  d}
\le 2 ( {2\a\log (m-n) \over m-n} )^{1/2}\atop 0< {\pi h\over \d}\le  ( {2\a\log
(m-n) \over m-n} )^{1/2}}     e^{-   {3n}  (  
{j\over d} -
 {h\over \d})^2    } \bigg\} . $$
 This establishes the lemma.\cqfd
 \medskip

     {\bf Case III}.  
Let 
 $${\bf \Psi}_{\bf 3}:= {1\over d\d}  \sum_{2 ( {2\a\log
(m-n) \over m-n} )^{1/2}\le {\pi j\over  d} \le \pi/2\atop 0< {\pi h\over \d}\le  ( {2\a\log
(m-n) \over m-n} )^{1/2}}   e^{ i\pi  (  \eta{j\over d}n +   \e{h\over \d}m) }  \cos^{m-n}{
\e   {\pi h\over
\d} } \cos^n{ \pi  ( \eta {j\over d} +  
\e{h\over \d})   } .$$

  \medskip \noi{\gem Lemma 2.4.}  {\it There exists   constants $C $ and $n_0$  depending on   $c$  such that for $n\ge n_0$
and 
$   n+ n_0 \le m\le 2n$}
$$   |{\bf \Psi}_{\bf 3}|  \le
C  
{1\over (m-n)^{  \a' }}({
\log m-n\over m-n})^{1/2}  . $$  
It is in this part that it is necessary to introduce the restriction $m\le 2n$.  
 \medskip
\noi{\it Proof.} 
 Here we have 
  
  $$\big( {2\a\log (m -n) \over  m -n }\big)^{1/2} \le \pi  (  {j\over d}  \pm    {h\over \d} )\le {\pi\over 2}+  \big( {2\a\log (m
 -n)\over m-n}\big)^{1/2} . $$ 
 Thus, as soon as $m-n$ is large enough, which is realized if $n_0$ is large enough,  we have   (as in case I) 
  the bound $|\cos \pi  (   {j\over d}  \pm    {h\over \d} )| \le \cos   \p_{m-n}   $. And we get  
 $|\cos \pi  (   {j\over d}  \pm    {h\over \d} )|^{ n} \le  \cos^{ n} \p_{m-n}    
\le e^{-2n\sin^2(\p_{m-n}/2)} $.  But    since $m\le 2n$ we have $n\ge m-n$, and so
$$\eqalign{ 2n \sin^2(\p_{m-n}/2) &=\big({n\over m-n}\big)2(m-n) \sin^2(\p_{m-n}/2)\ge 2(m-n) \sin^2(\p_{m-n}/2)\cr & = 2(m-n)
(\p_{m-n}/2)^2\tau_{m-n}^2\cr &
\ge     {\a'(m-n) \log m-n\over m-n} =  \a'  \log m-n 
    .\cr}$$

Therefore $\cos^{ n}
\p_{m-n} 
 \le (m-n)^{- \a'}$, and we have
 $$\eqalign{  {1\over d\d}  \sum_{2  
 ( {2\a\log (m -n) \over  m -n } )^{1/2}\le {\pi j\over  d}  \le \pi/2 \atop 0< {\pi
h\over \d}\le  ( {2\a\log (m-n) \over m-n} )^{1/2}}   
 |\cos { \pi  (  {j\over d} \pm {h\over \d})   } |^n& \le {1\over d\d(m-n)^{  \a' }} \sum_{2  
 ( {2\a\log (m -n) \over  m -n } )^{1/2}\le {\pi j\over  d}  \le \pi/2 \atop 0< {\pi
h\over \d}\le  ( {2\a\log (m-n) \over m-n} )^{1/2}} 1 \cr&  \le {1\over  \d(m-n)^{  \a' }} \sum_{  0< {\pi
h\over \d}\le  ( {2\a\log (m-n) \over m-n} )^{1/2}} 1
\cr&\le C
{1\over (m-n)^{  \a' }}({
\log m-n\over m-n})^{1/2} ,\cr} $$
as required.  \cqfd
  Now we estimate ${\bf r}$ defined in (2.8a), when $d$ is even. We have
$$|{\bf r} |\le {1\over d\d}   \sum_{1\le |h|<\d/2}    |\cos^{m-n}{     {\pi h\over
\d} } |.|\sin^n         
{\pi h\over \d}  | 
 $$   
--- If $\big( {2\a\log (m-n) \over m-n}\big)^{1/2}\le {\pi h\over \d} \le \pi/2 $, then
$$|{\bf r} |\le {1\over d\d}   \sum_{1\le |h|<\d/2}    |\cos^{m-n}{     {\pi h\over
\d} } | \le   d^{-1}   (m-n)^{-\a'}
 $$
--- If $0\le {\pi h\over \d}\le \big( {2\a\log (m-n) \over m-n}\big)^{1/2} $, then
$$|{\bf r} |\le {1\over d\d}   \sum_{1\le |h|<\d/2}   |\sin^n         
{\pi h\over \d}  | 
 \le   d^{-1} \big( {2\a\log (m-n) \over m-n}\big)^{n/2}  
. $$
 So that, for $n$ large enough   
$$ |{\bf r} |\le  d^{-1}   (m-n)^{-\a'}.\eqno(2.19)$$
\medskip\noi {\it Proof of Proposition 2.1.}   
Combining Lemmas 2.2, 2.3 and 2.4, finally gives in view of (2.8), estimate (1.3), and that  $\a >\a' >\max(3/2, 1/c)$: there exist
constants $C$ and $n_0$ such that for
 $n\ge n_0$, and 
   $      n+ n^{c} \le m\le 2n   $
$$\displaylines{\big|{\bf \Psi} \big( (d,n), (\d, m)\big)\big|\le  C \Bigg\{{1\over (m-n)^{  \a' }}+ n ({  \log (m-n)\over m-n})^{
1/c}\hfill\cr\hfill+
  {1\over d\d}\sum_{0< 
{\pi j\over  d}
\le 2 ( {2\a\log (m-n) \over m-n} )^{1/2}\atop 0< {\pi h\over \d}\le  ( {2\a\log
(m-n) \over m-n} )^{1/2}}     e^{-   {3n}  (  
{j\over d} -
 {h\over \d})^2    }  
+ 
{1\over (m-n)^{  \a' }}({
\log m-n\over m-n})^{1/2} \Bigg\}\cr \hfill \le  C \bigg\{  n ({  \log (m-n)\over m-n})^{ 1/c}+
{1\over d\d}  \sum_{0<
{\pi j\over  d}
\le 2 ( {2\a\log (m-n) \over m-n} )^{1/2}\atop 0< {\pi h\over \d}\le  ( {2\a\log
(m-n) \over m-n} )^{1/2}}     e^{-   {3n}  (  
{j\over d} -
 {h\over \d})^2    } 
 \bigg\}.   \quad\qq\  \cr}$$

\cqfd
{\bf Remarks.} --- The condition $m\ge n+n^c$ is in turn only used to   make  the bound in  Lemma 2.3 efficient.

\noi --- By construction, we have
the  trivial bound
$|{\bf
\Psi_2}|\le C { 
\log (m-n)\over m-n}
$. Combining it with Lemmas 2.2   and 2.4, we get another estimate: $|{\bf \Psi }|\le C {  \log (m-n)\over m-n} $  which is valid as soon
as
$m-n\ge n_0$,
$m\le 2n$,
$n_0$ sufficiently large.
\medskip  Turning to estimates involving the correlation ${\bf \Delta}$, we notice that $m\le 2n$ implies $m-n\le n$ and so 
$${\pi \over2\sqrt{2\a}   } ( {m-n\over \log (m-n)   } )^{1/2}\le  { \pi\over \sqrt{2\a}  }( { n \over \log n } )^{1/2}.$$
If $\max(d, \d) <  { \pi\over \sqrt{2\a}  }( { n \over \log n } )^{1/2} $ then by (1.16)  
$$ \big|\P\big\{ d | B_n\big\}-  {1\over d}
\big|\big|\P\big\{ \d | B_n\big\}-  {1\over \d}
\big|\le  Cn^{-2\a'}.  $$

Combining this with Proposition 2.1, shows in view of (2.11) \medskip\noi {\gem Proposition 2.5.} {\it There exist constants $n_0$ and  $C
$ such that for any
$  n\ge n_0$ and 
  $      n+ n^{c} \le m\le 2n   $
}
  $$\displaylines{ \sup_{ (d \vee \d) <  { \pi\over \sqrt{2\a}  }( { n \over \log n } )^{1/2}}\big|{\bf \Delta} \big( (d,n), (\d,
m)\big)\big| \le   C n  ({  \log (m-n)\over m-n})^{ 1/c} 
      . \cr}$$

\bigskip\centerline{\bf Case: ${\bf    m\ge 2n.} $}
\medskip\noi 
To treat this case, we have to proceed to little changes, but the method is very similar.  
We will
establish the following 
\medskip\noi {\gem Proposition 2.6.} {\it  There exist constants $C$ and $n_0$  depending on   $c$  such
that for
 $n\ge n_0$, and 
    $    m\ge  2n $  
$$  \big|{\bf \Psi} \big( (d,n), (\d, m)\big)\big| \le   C \bigg\{  
n({  \log n\over n})^{{1/ 2c}}({ \log (m-n)\over
m-n})^{1/2}+    {1\over d\d} \sum_{0<  {\pi j\over  d}
\le  \ ( {2\a\log n \over n} )^{1/2}\atop 0< {\pi h\over \d}\le  ( {2\a\log (m-n) \over m-n} )^{1/2}}    e^{-   3n (   {j\over d} -
 {h\over \d})^2   }     \bigg\}.    $$ 
Let $\d_0>\sqrt 2 $, then for all $n$ large and $m\ge 2n$}
$$\sup_{d<{\pi\over \sqrt{2\a}}  \sqrt{n\over \log n}  \atop \d<{\pi\over \d_0\sqrt{2\a}} \sqrt{m\over\log m}  }\big|{\bf
\Psi} \big( (d,n), (\d, m)\big)\big| \le  C    n({  \log n\over n})^{{1/ 2c}}({ \log (m-n)\over
m-n})^{1/2} $$

\medskip Before giving the proof, observe by the choice   of  $\a,\a'  $  that $\a >\a'> 
   \max\big(3/2,(1/ 2c)-1\big) $.  
 Next let $n_0 $   be sufficiently large  so that for $n\ge n_0$,
 $$\min \big( \tau_n, \tau_{m-n}\big)\ge
(\a^\prime/\a)^{1/2}  .\eqno(2.20)$$    
  We distinguish again between three cases: 
\smallskip \noi \qq\qq\qq\qq  --- (I) $\big( {2\a\log (m-n) \over m-n}\big)^{1/2}\le {\pi h\over \d} \le \pi/2 $.
\smallskip \noi \qq\qq\qq\qq  --- (II) $0< {\pi h\over \d}\le \big( {2\a\log (m-n) \over m-n}\big)^{1/2}$ and $0< 
{\pi j\over  d}
\le  2\big( {2\a\log n \over n}\big)^{1/2}$.
\smallskip \noi \qq\qq\qq\qq  --- (III) $0< {\pi h\over \d}\le \big( {2\a\log (m-n) \over m-n}\big)^{1/2}$ and $ 2\big( {2\a\log
n \over n}\big)^{1/2}\le {\pi j\over  d} \le \pi/2 $.
 
\medskip    {\bf Case I}. The sum
$${\bf \Psi}_{\bf 1}:= {1\over d\d}  \sum_{1\le j<{d\over 2}}\sum_{ ( {2\a\log (m-n) \over m-n} )^{1/2}\le {\pi h\over \d} \le \pi/2 }   e^{ i\pi 
( 
 {\eta j\over d}n +   {\e h\over \d}m) }  \cos^{m-n}{
\pi   {h\over
\d} } \cos^n{ \pi  (  {\eta j\over d} +  
{\e h\over \d})   }, $$
 has been already estimated in Lemma 2.2 and we recall that we have $ |{\bf \Psi}_{\bf 1}| \le (m-n)^{- \a'}$.

\medskip    {\bf Case II}.  Let
$${\bf \Psi}_{\bf 2}:={1\over d\d}  \sum_{0< 
{\pi j\over  d}
\le  2 ( {2\a\log n \over n} )^{1/2}}   \sum_{0< {\pi h\over \d}\le  ( {2\a\log (m-n) \over m-n} )^{1/2}}\!\!\!e^{ i\pi  (  \eta{j\over
d}n +
\e  {h\over \d}m) } 
\cos^{m-n}{ \pi   {h\over
\d} }\ \cos^n{ \pi  ( \eta {j\over d} + \e 
{h\over \d})   }     .  $$
  
We will establish
\medskip\noi {\gem Lemma 2.7.} {\it There exist    constants $C $ and $n_0$  depending on $c$, such that for $n\ge n_0$ and $m\ge 2n$} 
$$  |{\bf \Psi}_{\bf 2}
|\le 
       C\Big(n({  \log n\over n})^{{1/ 2c}}({ \log (m-n)\over m-n})^{1/2}+ {1\over d\d} \sum_{0< 
{\pi j\over  d}
\le  2 ( {2\a\log n \over n} )^{1/2}\atop 0< {\pi h\over \d}\le  ( {2\a\log (m-n) \over m-n} )^{1/2}}    e^{-   3n  (   {j\over d} -
 {h\over \d})^2    }     \Big) .   $$

\medskip\noi {\it Proof.} We proceed as before except that we select $p$ so that 
 $ c> {1\over 2p +1}$,  
and will compare the sum ${\bf \Psi}_{\bf 2}$
  with the sum
 $$ {\bf \Psi}_{\bf 2}':={1\over d\d}  \sum_{0<
{\pi j\over  d}
\le  2 ( {2\a\log n \over n} )^{1/2}}   \sum_{0< {\pi h\over \d}\le  ( {2\a\log (m-n) \over m-n} )^{1/2}}  e^{ i\pi  (  \eta{j\over d}n + \e 
{h\over \d}m) } 
\cos^{m-n}{ \pi   {h\over
\d} } E^n(\pi  (  {\eta j\over d} +  
{\e h\over \d}))    . \eqno(2.21)$$
  Again we   observe that 
 $\big|\cos^n x -e^{ nQ_p(x)}\big|\le n\big|\log\cos  x - Q_p(x)\big|\le C_pn|x|^{2p}$,  
for   $|x|\le x_p$, and that for $n$ large enough we have  $\pi | {\eta j\over d} +
 {\e h\over \d}|\le x_p$.    Thus
$$  \big|\cos^n{ \pi  (  {\eta j\over d} +  
{\e h\over \d})   }- E^n(\pi  (  {\eta j\over d} +  
{\e h\over \d})) \big| \le  C_pn\big[ ({h\over \d})^{2p} +    ({j\over
d})^{2p} \big].  $$
  And so
$$\bigg|{1\over d\d}  \sum_{0<
{\pi j\over  d}
\le   2 ( {2\a\log n \over n} )^{1/2}\atop 0< {\pi h\over \d}\le  ( {2\a\log (m-n) \over m-n} )^{1/2}} e^{ i\pi  (  \eta{j\over d}n +   \e{h\over \d}m) } 
\cos^{m-n}{
\e   {\pi h\over
\d} }\Big\{\cos^n{ \pi  ( \eta {j\over d} +  
\e{h\over \d})   }-E^n(\pi  (  {\eta j\over d} +  
{\e h\over \d})) \Big\}\bigg|$$
$$\eqalign{& \le C {n\over d\d}   \sum_{0<
{\pi j\over  d}
\le 2 ( {2\a\log n \over n} )^{1/2}\atop 0< {\pi h\over \d}\le  ( {2\a\log (m-n) \over m-n} )^{1/2}}
      \big[ ({h\over \d})^{2p} +    ({j\over
d})^{2p} \big] .\cr} $$ 
 Therefore
$$ \big|{\bf \Psi}_{\bf 2}-{\bf \Psi}'_{\bf 2}\big|\le  Cn \max \bigg(({ \log n\over n})^{1/2}({  \log (m-n)\over
m-n})^{( p+1 /2)},({ 
\log n\over n})^{( p+1 /2)}({
\log (m-n)\over m-n})^{1/2}\bigg).\eqno(2.22)$$
 
\noi As we have 
$$  ({ \log n\over n})^{1/2}({  \log (m-n)\over m-n})^{( p+1 /2)}\ge ({  \log n\over n})^{( p+1 /2)}({ \log (m-n)\over m-n})^{1/2}\qquad
\Leftrightarrow\qquad m-n\le n.\eqno(2.23)$$
Indeed, the   inequality on the left is equivalently rewritten  as $({  \log (m-n)\over m-n})^p\ge ({ \log n\over
n})^p$.
 Further, since  
 $E(x)  \le e^{-   x ^2 /3}  $ for $|x|\le x_p$, we may again simply bound $|{\bf \Psi}_{\bf 2}'|$    
as follows: for $n$ large enough
 
$$\eqalign{  |{\bf \Psi}_{\bf 2}'|\le {1\over d\d}  \sum_{0< 
{\pi j\over  d}
\le 2 ({2\a\log n \over n} )^{1/2}\atop 0< {\pi h\over \d}\le  ( {2\a\log (m-n) \over m-n} )^{1/2}}   e^{-  {3n} (  {\eta j\over d} +
 {\e h\over \d})^2 }  .
 \cr}\eqno(2.24)$$  
   It is clear that   the  sum appearing in the righthand side is bounded by   
$$   {1\over d\d}   \sum_{0<
{\pi j\over  d}
\le 2 ({2\a\log n \over n} )^{1/2}\atop 0< {\pi h\over \d}\le  ( {2\a\log (m-n) \over m-n} )^{1/2}}     e^{-   {3n}  (  
{j\over d} -
 {h\over \d})^2    }     \eqno(2.25)$$
    With  estimates (2.22), (2.23), (2.24) and (2.25), we therefore get if $ m\ge 2n $ 
$$  |{\bf \Psi}_{\bf 2}
|\le 
     \displaystyle{ C\Big(n({  \log n\over n})^{( p+1 /2)}({ \log (m-n)\over m-n})^{1/2}+ {1\over d\d}  \sum_{0<
{\pi j\over  d}
\le 2 ({2\a\log n \over n} )^{1/2}\atop 0< {\pi h\over \d}\le  ( {2\a\log (m-n) \over m-n} )^{1/2}}    e^{-   {3n}  (   {j\over d} -
 {h\over \d})^2    }     \Big)}    . 
 $$
As $ c> {1\over 2p +1}$, it follows that 
$$  |{\bf \Psi}_{\bf 2}
|\le 
       C\Big(n({  \log n\over n})^{{1/ 2c}}({ \log (m-n)\over m-n})^{1/2}+ {1\over d\d} \sum_{0< 
{\pi j\over  d}
\le  2 ( {2\a\log n \over n} )^{1/2}\atop 0< {\pi h\over \d}\le  ( {2\a\log (m-n) \over m-n} )^{1/2}}    e^{-   3n  (   {j\over d} -
 {h\over \d})^2    }     \Big) .   $$
 This establishes the lemma.\cqfd
 \medskip

     {\bf Case III}.   
Let 
 $${\bf \Psi}_{\bf 3}:= {1\over d\d}  \sum_{  2( {2\a\log
n \over n} )^{1/2}\le {\pi j\over  d} \le \pi/2 }\sum_{0< {\pi h\over \d}\le  ( {2\a\log (m-n) \over m-n} )^{1/2}}  \!\!\!\!\!\!\!\!\! e^{
i\pi  ( 
\eta{j\over d}n +   \e{h\over \d}m) }  \cos^{m-n}{
\e   {\pi h\over
\d} } \cos^n{ \pi  ( \eta {j\over d} +  
\e{h\over \d})   } .$$
  \medskip \noi{\gem Lemma 2.8.} {\it }
$$   |{\bf \Psi}_{\bf 3}|  \le
{1\over n^{  \a' }}({
\log m-n\over m-n})^{1/2} . $$  
 
 \noi{\it Proof.}  
 Here we have since $m-n\ge n$
  
  $$ \eqalign{    {\pi\over 2}+  \big( {2\a\log (m
 -n)\over m-n}\big)^{1/2} &\ge \pi  (  {j\over d}  \pm    {h\over \d} )\cr &\ge  2( {2\a\log
n \over n} )^{1/2}-\big( {2\a\log (m -n) \over  m -n }\big)^{1/2}\ge  ( {2\a\log
n \over n} )^{1/2} .\cr} $$ 
 Thus, as soon as   $n$ is large enough,  we have   
  the bound $|\cos \pi  (   {j\over d}  \pm    {h\over \d} )| \le \cos   \p_{ n}   $. And we get  
 $|\cos \pi  (   {j\over d}  \pm    {h\over \d} )|^{ n} \le  \cos^{ n} \p_{ n}    
=e^{-2n\sin^2(\p_{ n}/2)} $.  But    
$$\eqalign{ 2n \sin^2(\p_{ n}/2) & \ge 2n\sin^2(\p_{ n}/2)  = 2n
(\p_{n}/2)^2\tau_{n}^2 
\ge     {\a'n\log n\over n} =  \a'  \log n 
    .\cr}$$ 

\noi Therefore $\cos^{ n}
\p_{n} 
 \le n^{- \a'}$, and so have
 $$\eqalign{  {1\over d\d}  \sum_{2  
 ( {2\a\log n \over n } )^{1/2}\le {\pi j\over  d}  \le \pi/2 \atop 0< {\pi
h\over \d}\le  ( {2\a\log (m-n) \over m-n} )^{1/2}}   
 |\cos { \pi  (  {j\over d} \pm {h\over \d})   } |^n& \le {1\over d\d n^{  \a' }} \sum_{2  
 ( {2\a\log n \over  n } )^{1/2}\le {\pi j\over  d}  \le \pi/2 \atop 0< {\pi
h\over \d}\le  ( {2\a\log (m-n) \over m-n} )^{1/2}} 1 \cr&  \le {1\over  \d n^{  \a' }} \sum_{  0< {\pi
h\over \d}\le  ( {2\a\log (m-n) \over m-n} )^{1/2}} 1
\cr&\le C
{1\over n^{  \a' }}({
\log m-n\over m-n})^{1/2} ,\cr} $$
as claimed.  \cqfd

  Finally ${\bf r}$ is estimated in exactly the same way in this case too, and we have that estimate (2.19) holds again. 
 
\medskip\noi {\it Proof of Proposition 2.6.}    Combining
the previous estimates finally show in view of (2.8), estimate (1.3) and since  $\a >\a' >\max(3/2, (1/c)-1)$: there exist constants $C$
and $n_0$ such that for
 $n\ge n_0$, and 
    $    m\ge  2n $ 
$$\displaylines{ \big|{\bf \Psi} \big( (d,n), (\d, m)\big)\big| \le  C \Bigg\{ {1\over (m-n)^{  \a'} } +   
n({  \log n\over n})^{{1/ 2c}}({ \log (m-n)\over
m-n})^{1/2}+  \hfill\cr\hfill +{1\over d\d} \sum_{0<  {\pi j\over  d}
\le  \ ( {2\a\log n \over n} )^{1/2}\atop 0< {\pi h\over \d}\le  ( {2\a\log (m-n) \over m-n} )^{1/2}}    e^{-   3n (   {j\over d} -
 {h\over \d})^2   }     +{1\over n^{  \a' }}({
\log m-n\over m-n})^{1/2} \Bigg\}\cr \hfill  \le  C \bigg\{  
n({  \log n\over n})^{{1/ 2c}}({ \log (m-n)\over
m-n})^{1/2}+    {1\over d\d} \sum_{0<  {\pi j\over  d}
\le  \ ( {2\a\log n \over n} )^{1/2}\atop 0< {\pi h\over \d}\le  ( {2\a\log (m-n) \over m-n} )^{1/2}}    e^{-   3n (   {j\over d} -
 {h\over \d})^2   }     \bigg\}. \quad\ (2.26)\cr} $$
And produces the wished inequality. Now if $d<{\pi\over \sqrt{2\a}}\big({n\over \log n}\big)^{1/2}$, 
$\d<{\pi\over \d_0\sqrt{2\a}}\big({m\over\log m}\big)^{1/2}$,   as $m\le 2(m-n)$ we have
$$\d<{\pi\over \d_0\sqrt{2\a}}\big({2(m-n)\over\log 2(m-n)}\big)^{1/2} <{\pi \over  \sqrt{2 \a}}\big({ (m-n)\over\log
 (m-n)}\big)^{1/2},$$
if $   \d_0 >\sqrt 2$. In which case, the exponential sum appearing in the righthand side of (2.26)   no longer contributes. Hence
$$\sup_{d<{\pi\over \sqrt{2\a}}  \sqrt{n\over \log n}  \atop \d<{\pi\over \d_0\sqrt{2\a}} \sqrt{m\over\log m}  }\big|{\bf
\Psi} \big( (d,n), (\d, m)\big)\big| \le  C    n({  \log n\over n})^{{1/ 2c}}({ \log (m-n)\over
m-n})^{1/2}.\eqno(2.27) $$
\cqfd
{\bf Remark.}   --- We have
the  trivial bound
$|{\bf
\Psi_2}|\le C ( { \log n \over n} )^{1/2} ( { \log (m-n) \over m-n} )^{1/2}
$. By combining   with Lemmas 2.6   and 2.8  we also get: $|{\bf \Psi }|\le  ( { \log n \over n} )^{1/2} ( { \log (m-n)
\over m-n} )^{1/2} $,  which is valid for
$m \ge 2n\ge n_0$,
 $n_0$ sufficiently large.
\medskip  Turning to estimates involving the correlation ${\bf
\Delta}$, we have in view of  
  (1.16)
$$ \big|\P\big\{ d | B_n\big\}-  {1\over d}
\big|\big|\P\big\{ \d | B_m\big\}-  {1\over \d}
\big|\le  Cn^{- \a'}m^{- \a'}.  $$
With (2.26)  this now implies
$$\eqalign{\sup_{d<{\pi\over \sqrt{2\a}}  \sqrt{n\over \log n}  \atop \d<{\pi\over \d_0\sqrt{2\a}} \sqrt{m\over\log m}  }\big|{\bf \Delta}
\big( (d,n), (\d, m)\big)\big|&\le C  
n({  \log n\over n})^{{1/ 2c}}({ \log (m-n)\over
m-n})^{1/2}+{C\over(nm)^{\a'} } 
\cr  &\le C  
n({  \log n\over n})^{{1/ 2c}}({ \log (m-n)\over
m-n})^{1/2}.\cr}$$
\medskip\noi {\gem Proposition 2.9.} {\it There exist two constants $C$ and $n_0$  depending on $c$, such that for any $  n\ge
n_0$  
   and $    m\ge  2n $}
 $$\sup_{d<{\pi\over \sqrt{2\a}}  \sqrt{n\over \log n}  \atop \d<{\pi\over \d_0\sqrt{2\a}} \sqrt{m\over\log m}  }\big|{\bf \Delta}
\big( (d,n), (\d, m)\big)\big| \le C  
n({  \log n\over n})^{{1/ 2c}}({ \log (m-n)\over
m-n})^{1/2}.$$

\medskip\noi\centerline {\gem 2.3. The strongly dependent case }

   \medskip\noi  
The main object of this section will consist of establishing the following delicate
estimate. 

\medskip\noi {\gem Proposition 2.10.} {\it For any $\e>0$, there exists a constant $C_\e$ depending on $\e$ only, such that 
$$\Big|\P \{ D|B_nB_m \}- {1\over D  2^{ m-n}} \sum_{k=0}^{m-n} { C^k_{m-n}     } \rho_k(D) \Big|\le C_\e \left({     D  ^{1 + \e}\over  
n }\right)^{1/2 },  $$
where 
$$\rho_k(D)=\cases{\displaystyle{\ \ \ \ \ \ \, \prod_{p|D} p^{\lfloor {v_p(D)\over 2}\rfloor }} \qq &\hbox{ if $k=0$,}\cr
\displaystyle{ \prod_{ v_p(k)<{v_p(D)/
2}}(2p^{v_p(k)})\cdot\prod_{   v_p(k)\ge {v_p(D)/
2} }  p^{\lfloor{v_p(D)\over
2}\rfloor}} \qq &\hbox{ if $k\ge 1$.}\cr} $$Further, for any positive integers $n,m$ and $D$,  
$$  \P \{ D|B_nB_m \}    \le   {1\over 2^{  m-n }\sqrt D
}+    {2^{\o(D)}( m-n)  \over D } + C_\e
 {     D  ^{{1/ 2} + \e}\over   \sqrt n }  .  
   $$
And for any $\e>0$, there exists a constant $C_\e$ depending on $\e$ only, such that} 
$$\P \{ 
d|B_n\, ,\, \d|B_m \}\le  C   (m-n){ 2^{\o(d\d) } \over d\d
}+ C_\e  {     (d\d)  ^{(1 + \e)/2}\over  \sqrt n } .$$

\medskip
The proof of Proposition 2.10 is based on several intermediate results.   We   begin with
computing the characteristic function of $B_nB_m$. Plainly, writing that $B_nB_m= B_n^2 +B_n(B_m-B_n)$ and using independence
$$\E e^{i\upsilon B_nB_m}={1\over 2^n2^{m-n}}\sum_{k=0}^{m-n}C^k_{m-n}\sum_{h=0}^n   C^h_n e^{i\upsilon (h^2+kh)}.$$
And so 
$$\P \{ D|B_nB_m \}=\E\Big({1\over D} \sum_{j=0}^{D-1} e^{2i\pi  {j \over D}B_nB_m}\Big)= {1\over 2^n2^{m-n}}\sum_{k=0}^{m-n}C^k_{m-n}{1\over
D}\sum_{j=0}^{D-1}\sum_{h=0}^n   C^h_n e^{2i\pi{j\over D}  (h^2+kh)}. \eqno(2.28)$$
 We preliminarily evaluate the $({\bf C}, 1)$ sums $$S_L:={1\over L}\sum_{h=0}^L  e^{2i\pi{j\over D} 
(h^2+kh)},$$
and   will next compare $S_n$ to the Euler $({\bf E}, 1)$ sum    $\sum_{h=0}^n   2^{-n}C^h_n e^{2i\pi{j\over D} (h^2+kh)}$. We write
$L=ND+m$ with $N \ge 0$ and 
$0\le m< D$. Then
$$S_L=  {1\over L}\left(\sum_{X=0}^{N-1} \sum_{x=XD+1}^{(X+1) D } + \sum_{x=ND+1}^{ND+m}\right)e^{2i\pi{j\over D} 
(x^2+kx)}$$
As $$ \sum_{x=XD+1}^{(X+1) D }e^{2i\pi{j\over D} 
(x^2+kx)}=\sum_{y= 1}^{D}e^{2i\pi{j\over D} 
(y^2+ky)},$$
we have
$$\eqalign{S&={N\over L}\sum_{y= 1}^{D}e^{2i\pi{j\over D} 
(y^2+ky)}+ {1\over L}\sum_{y= 1}^{m}e^{2i\pi{j\over D} 
(y^2+ky)}
\cr&=N\big({1\over L}-{1\over ND}\big)\sum_{y= 1}^{D}e^{2i\pi{j\over D} 
(y^2+ky)}+ {1\over D}\sum_{y= 1}^{D}e^{2i\pi{j\over D} 
(y^2+ky)}+ {1\over L}\sum_{y= 1}^{m}e^{2i\pi{j\over D} 
(y^2+ky)}\cr}$$
Therefore, for   $j=1,\ldots, D$ 
$$\eqalign{\Big|S_L-{1\over D}\sum_{y= 1}^{D}e^{2i\pi{j\over D} 
(y^2+ky)} \Big|&\le  \big( {|ND-L|\over L D}\big)\big|\sum_{y= 1}^{D}e^{2i\pi{j\over D} 
(y^2+ky)} \big|  + {1\over L}\big|\sum_{y= 1}^{m}e^{2i\pi{j\over D} 
(y^2+ky)}\big|
\cr &\le    {1\over  L} \Big(\big|\sum_{y= 1}^{D}e^{2i\pi{j\over D} 
(y^2+ky)} \big|  +  \big|\sum_{y= 1}^{m}e^{2i\pi{j\over D} 
(y^2+ky)}\big|\Big)\cr} \eqno(2.29)$$

\noi Note that if $t=(j,D)>1$,   estimate (2.29) can be improved by using the same arguments. Let $j=t j' $, $D=tD'$. Then $ {1\over
D}\sum_{y= 1}^{D}e^{2i\pi{j\over D}  (y^2+ky)} ={1\over D'}\sum_{y= 1}^{D'}e^{2i\pi{j\over D'} 
(y^2+ky)}$,   and   writing $L$ under the form $L= N'D' +m'$ with $0\le m'<
D'$, we get similarly    
$$\eqalign{\Big|S_L-{1\over D }\sum_{y= 1}^{D }e^{2i\pi{j \over D } 
(y^2+ky)} \Big|&\le {1\over L}\Big(\big|\sum_{y= 1}^{D'}e^{2i\pi{j'\over D'} 
(y^2+ky)} \big|  +  \big|\sum_{y= 1}^{m}e^{2i\pi{j'\over D'} 
(y^2+ky)}\big|\Big)\cr &\le {2\over L}\max_{m'=1}^{D'} \big|\sum_{y=1}^{m'} e^{2i\pi{j'\over D'} 
(y^2+ky)} 
\big|.
\cr} \eqno(2.30) $$
When $L=n$,   this shows that 
$${1\over n}\sum_{h=0}^n  e^{2i\pi{j\over D} 
(h^2+kh)}={1\over D}\sum_{y= 1}^{D}e^{2i\pi{j\over D} 
(y^2+ky)} +{\cal O}\big( {1\over  n}\big). \eqno(2.31)$$
 \medskip Although     ({\bf
E},1) does not includes ({\bf C},1) (see  [H]  Chap.8) in general, the latter estimate will imply this, thanks to the result
below.
  \medskip
\noi{\gem  Lemma 2.11.}  ([H] Theorem 149, p.213) {\it Let $A=\{ A_n, n\ge 1\}$ be a sequence of reals such that
 $$C^1_n(A)= {A_1+\ldots +A_n\over n}= a + o(n^{-1/2}). $$
Then $A$ is summable $(\E, q)$ for any positive $q$.}
\medskip
\noi The conclusion of the lemma is wrong when replacing 
$o$ by 
${\cal O}$ (see also [H]).
\medskip Let $$\rho_k(D)=\#\big\{1\le y\le D: \  D| y^2+ky \big\} ,\qq \quad k=0,1,\ldots, m-n .  $$ 
 
  \medskip
\noi{\gem  Corollary 2.12.}   {\it We have for each $k$} 
$$  \lim_{n\to\infty}  {1\over D}\sum_{j=0}^{D-1}  \sum_{h=0}^n   2^{-n}C^h_n e^{2i\pi{j\over D} (h^2+kh)}={\rho_k(D)\over D }
  .  $$  \medskip
\noi{\it Proof.} In view of (2.29), the assumption of   lemma 2.15 is fulfilled.
 Thus, by considering separately imaginary and real parts,   the lemma applied with
$q=1$    implies for  $j=1,\ldots, D$ that
$$     \lim_{n\to\infty}  \sum_{h=0}^n   2^{-n}C^h_n e^{2i\pi{j\over D} (h^2+kh)}=  {1\over
D}\sum_{y= 1}^{D}e^{2i\pi{j\over D}  (y^2+ky)}  .  $$
Henceforth   
$$  \lim_{n\to\infty}   {1\over D}\sum_{j=0}^{D-1}  \sum_{h=0}^n   2^{-n}C^h_n e^{2i\pi{j\over D} (h^2+kh)}={1\over D^2}
\sum_{j=0}^{D-1}\sum_{y= 1}^{D}e^{2i\pi{j\over D}  (y^2+ky)}  . $$ It remains to observe that 
$$ {1\over  D^2}\sum_{j=0}^{D-1} \sum_{y= 1}^{D}e^{2i\pi{j\over D} 
(y^2+ky)}={\#\big\{1\le y\le D: D|y^2+ky\big\}\over D }
={\rho_k(D)\over D } . \eqno(2.32)$$
\cqfd
   
We shall now give a uniform estimate of the speed of convergence in the above limit. For, we use the simple bound of the difference
between  the sum
({\bf E},1)  by the Euler method  and the sum ({\bf C},1) by the C\'esaro method. By linearity, it is enough to get a bound for the sum
({\bf E},1) alone. Let $\{a_k, k\ge
0\}$ be an arbitrary sequence of reals and put $A_\ell= \sum_{k=0}^\ell a_k$,  $\ell\ge 0$.
\medskip Then there exists an   absolute constant $C$ such that for every positive integer $n$
$$\big| \sum_{h=0}^n   2^{-n}C^h_n a_h\big|\le {C\over \sqrt n}\max_{\ell=0}^n\big|A_\ell\big|, \eqno(2.33)$$
  This is easily seen by using Abel summation:   put
$$E_n=\sum_{h=0}^n   2^{-n}C^h_n a_h=\sum_{h=0}^n   v_h a_h, $$
where  $v_h=v_h(n)=2^{-n}C^h_n$, $h=0,\ldots, n$. According to Theorem 138(1)
p.201 in  [H] (see also p. 214), the supremum is reached at the value
$$\nu = [{n+1\over 2}]\qq{\rm and}\qq v_\nu\le {C\over \sqrt n} , \eqno(2.34)$$
 where $C$ is an absolute constant. If ${n+1\over 2}$ is integer, then $v_{\nu-1}$ and $v_\nu$  are equal. Besides, $v_k$ decreases on
either side of $k=\nu$.  Splitting the sum $E_n $ into the two subsums  
 $ E^1_n = \sum_{k=0}^\nu v_ka_k, \ E^2_n =\sum_{k=\nu+1}^n v_ka_k   $, and since   $a_\ell= A_\ell-A_{\ell -1}$,
$a_0=A_0$, we get in the one hand 
$$ E^1_n  =v_0a_0+ \sum_{k=1}^\nu v_k(A_k-A_{k-1})     = -\big\{A_0(v_1-v_0)+ A_1(v_2-v_1)+\ldots+ A_{\nu-1}(v_\nu-v_{\nu-1})\big\} +v_\nu
A_\nu . $$
Thus 
$$ |E^1_n | \le v_\nu| A_\nu|+\max_{m=0}^{\nu-1}|A_m| \sum_{m=1}^\nu (v_m-v_{m-1}) \le 2v_\nu\max_{m=0}^{\nu }|A_m|. $$ 
   And in the other 
 $$ E^2_n  = \sum_{k=\nu+1}^n v_k(A_k-A_{k-1})   =-v_{\nu+1} A_{\nu }
+A_{\nu+1}(v_{\nu+1}-v_{\nu+2})+ \ldots+ A_{n-1}(v_{n-1}-v_{n})+v_nA_n.  $$ 
Hence
$$ |E^2_n | \le \max_{m=\nu}^{n}|A_m|\big\{v_{\nu+1}  +\sum_{m=\nu+1}^{n-1} (v_m-v_{m+1})
+v_n\big\}  \le 2v_{\nu+1}\max_{m=\nu}^{n}|A_m|   .  $$ 
 We have thus established  (2.33). 
\medskip Now let $0\le j\le D-1$. Let also $0\le k\le m-n$. We apply (2.33) with the choice 
$$a_h=  e^{2i\pi{j\over D} (h^2+kh)}- {1\over D}\sum_{y= 1}^{D}e^{2i\pi{j\over D} 
(y^2+ky)} , \qq 0\le h\le n.$$
If $j=0$, then $a_h\equiv 0$ and there is nothing to prove. If $0< j\le D-1$, we find in view of (2.30) that
 $$\eqalign{\qq & \Big|\sum_{h=0}^n   2^{-n}C^h_n \Big(e^{2i\pi{j\over D} (h^2+kh)}- {1\over D}\sum_{y= 1}^{D}e^{2i\pi{j\over D} 
(y^2+ky)} \Big)\Big|\cr  &= \Big|\sum_{h=0}^n   2^{-n}C^h_n e^{2i\pi{j\over D} (h^2+kh)}- {1\over D}\sum_{y= 1}^{D}e^{2i\pi{j\over D} 
(y^2+ky)} \Big|
\cr &\le {C\over \sqrt n}\max_{\ell=0}^n\big|\sum_{h=0}^\ell\Big(e^{2i\pi{j\over D} (h^2+kh)}- {1\over D}\sum_{y= 1}^{D}e^{2i\pi{j\over D} 
(y^2+ky)}\Big) \big|
\cr &=  {C\over \sqrt n}\max_{\ell=0}^n\big|\sum_{h=0}^\ell e^{2i\pi{j\over D} (h^2+kh)}- {\ell\over D}\sum_{y= 1}^{D}e^{2i\pi{j\over D} 
(y^2+ky)}  \big|\le  2{C\over \sqrt n}\max_{m'=1}^{D'} \big|\sum_{y=1}^{m'} e^{2i\pi{j'\over D'} 
(y^2+ky)} 
\big|,\cr}\eqno(2.35)$$
   where $C$ is the same absolute constant as in (2.33) and   the  used notation arises from (2.30): if $t=(j,D)>1$, then
$j=t j'
$,
$D=tD'$, $n= N'D' +m'$ with $1\le m'\le
D'$.
\medskip  
Now we need the following lemma: 
\medskip\noi{\gem Lemma 2.13.}   {\it Let  $\a$ be a real number and    $a,q$ be positive integers such that $(a,q)=1$ and
$|\a-a/q|<1/q^2$. Then, for any positive integer  $M$}
$$\big|\sum_{x=1}^M e^{2i\pi \a x^2 } \big|^2\le \sum_{u=1-M}^{M-1}\Big|\sum_{y=\max(1-u,1) }^{\min(M,M-u)}e^{4i\pi \a uy }\Big|\le 49
\Big\{{M^2\over  q}+(M\log q) + q\log q \Big\}.
$$ 
This follows from Lemma 4 p.128 in [S] and its proof. The last inequality is precisely what is established in the proof. If
$T=\sum_{x=1}^M e^{2i\pi
\a (x^2 +kx)}$, we have 
$$T^2= \sum_{x=1}^M\sum_{y=1}^M e^{2i\pi \a (x^2-y^2 +k(x-y))}= \sum_{y=1}^M\sum_{x=1}^M e^{2i\pi \a (x -y)(x+y+  k )}=
\sum_{y=1}^M\sum_{u=1-y}^{M-y} e^{2i\pi \a u(u+2y+  k )}$$
$$=\sum_{u=1-M}^{M-1}\sum_{y=\max(1-u,1) }^{\min(M,M-u)}  e^{2i\pi \a u(u+2y+  k )}\le \sum_{u=1-M}^{M-1}\Big|\sum_{y=\max(1-u,1)
}^{\min(M,M-u)}  e^{4i\pi \a   uy }\Big|$$
So that in turn
$$\sup_{k\ge 0}\Big|\sum_{x=1}^M\sum_{y=1}^M e^{2i\pi \a (x^2-y^2 +k(x-y))}\Big|^2\le 49
\Big\{{M^2\over  q}+(M\log q) + q\log q \Big\},$$
or
$$\sup_{k\ge 0}\Big|\sum_{x=1}^M\sum_{y=1}^M e^{2i\pi \a (x^2-y^2 +k(x-y))}\Big| \le 7
\Big\{{M \over  \sqrt q}+\sqrt{ M\log q} + \sqrt{q\log q} \Big\}.$$

 Thus
$$\eqalign{\sup_{k\ge 0}\big|\sum_{y=1}^{m} e^{2i\pi{j'\over D'} 
(y^2+ky)} 
\big|&\le 7\Big\{{m \over  \sqrt{ D'}}+\sqrt{ m\log D'} + \sqrt{D'\log D'} \Big\}
 \cr (m \le
D')\quad & \le  7  \Big\{ \sqrt{ D'}+2\sqrt{ D'\log D'} \Big\}
\cr   & \le C_\e  (D')^{1/2+ \e}\le C_\e  (D )^{1/2+ \e}. \cr} $$

Inserting this estimate  into (2.35)   leads to 
$$\sup_{ 0\le j< D }\sup_{0\le k\le m-n}\Big|\sum_{h=0}^n   2^{-n}C^h_n e^{2i\pi{j\over D} (h^2+kh)}- {1\over D}\sum_{y=
1}^{D}e^{2i\pi{j\over D}  (y^2+ky)} \Big|\le C_\e \left({     D  ^{1 + \e}\over   n }\right)^{1/2 }. \eqno(2.36)$$
Thereby   in view of  (2.32), (2.36) 
 $$  \sup_{0\le k\le m-n} \bigg|{1\over D}\sum_{   j=0}^{D-1}\sum_{h=0}^n   2^{-n}C^h_n e^{2i\pi{j\over D} (h^2+kh)}
   -{\rho_k(D)\over D } \bigg|  \le C_\e \left({     D  ^{1 + \e}\over   n }\right)^{1/2 }\!\!
.  \eqno(2.37)  $$ If now
we combine (2.37) with (2.28), we obtain
\medskip\noi {\gem Proposition 2.14.}  $$ \Big|\P \{ D|B_nB_m \}-   \sum_{k=0}^{m-n} { C^k_{m-n} \over 2^{ (m-n)}  } {\rho_k(D)\over D } 
\Big|
\le C_\e
\left({     D  ^{1 + \e}\over   n }\right)^{1/2 }. 
   $$
 \cqfd
\medskip 
{\gem Remarks.} 
\smallskip 
\item 1. It is possible ([S2]) to replace the error term $D^\e$ in Proposition 2.14 by a $(\log D)^2$ factor.
\item 2. One can bound    the difference between ${\bf \Delta}
\big( (d,n), (\d, m)\big)$ and ${\bf \Delta} \big( (d,n), (\d, m')$  once $m,m'$ are not too close to $n$. Indeed,
one can prove that there exists $n_0 $, such that if $m'\ge m \ge n+n_0$, then
 $$\Big|{\bf \Delta} \big( (d,n), (\d, m)\big) -{\bf \Delta} \big( (d,n), (\d, m')\big)\Big|\le  C\big({\log (m-n)\over
m-n}\big)^{1/2}.\eqno(2.38)$$ 
  In view of Lemma 2.3, ${\bf \Delta} \big( (d,n) , (\d, m)\big) $ equals to
$$ {1\over d\d}\sum_{1\le |h|<\d/2\atop1\le |j|<d/2}    e^{ i\pi  (  {j\over d}n +   {h\over \d}m) }\cos^{m-n}{ \pi  
{h\over
\d} }  \Big\{  \cos^n{ \pi  (  {j\over d} +  
{h\over \d})   }  -    \cos^{n} {
\pi j\over d}  \cos^{n} {
\pi h\over \d}\Big\}  $$  
So that for $n\le m\le m'$,
$$\displaylines{{\bf \Delta} \big( (d,n), (\d, m)\big) -{\bf \Delta} \big( (d,n), (\d, m')\big)=  {1\over d\d}\!\sum_{1\le
|h|<\d/2\atop1\le |j|<d/2}    \Big\{\cos^n{ \pi  (  {j\over d} +   {h\over \d})   }-\cos^{n} {
\pi j\over d}     \cos^{n} {
\pi h\over \d}\Big\}\hfill\cr\hfill\times \ e^{ i\pi     {j\over d}n   }\bigg[e^{ i\pi     {h\over \d}m  }  \cos^{m-n}{ \pi  
{h\over
\d} }-e^{ i\pi     {h\over \d}m'  }  \cos^{m'-n}{ \pi   {h\over
\d} }\bigg].\cr}$$
 Then   $$\Big|{\bf \Delta} \big( (d,n), (\d, m)\big) -{\bf \Delta} \big( (d,n), (\d, m')\big)\Big|\le   {2\over d\d}\!\!\sum_{1\le
|h|<\d/2\atop1\le |j|<d/2}   \Big|  \cos^{m-n}{ \pi   {h\over
\d} }-e^{ i\pi     {h\over \d}(m'-m)  }  \cos^{m'-n}{ \pi   {h\over
\d} }\Big|.\eqno(2.39)$$ 
 
In view of (2.13) and Case I of the proof of Theorem 2.4, if ${\pi h\over \d} \in I'_{m-n}$, then 
    $|\cos {\pi h\over \d}|\le \cos \p_{m-n}$. And so, for some $n_0 $   sufficiently large, and 
  $m-n\ge n_0$, 
  $$|\cos {\pi h\over \d}|^{m'-n} \le|\cos {\pi h\over \d}|^{m-n} \le  \cos^{m-n} \p_{m-n}    
=e^{-2(m-n)\sin^2(\p_{m-n}/2)}\le (m-n)^{-\b'}  .$$ 
    
 Hence,  
$$\displaylines{ {2\over d\d}  \sum_{{1\le |j|<{d\over 2}\,,\,  1\le |h|<{\d\over 2}}\atop {\pi h\over \d} \in
I'_{m-n}}   \Big|  \cos^{m-n}{ \pi   {h\over
\d} }- e^{ i\pi     {h\over \d}(m'-m)  }   \cos^{m'-n}{ \pi   {h\over
\d} }\Big|\hfill\cr\hfill\le {2\over d\d} \sum_{{1\le |j|<{d\over 2}\,,\,  1\le |h|<{\d\over 2}}\atop {\pi h\over \d} \in
I'_{m-n}}   \big(  |  \cos^{m-n}{ \pi   {h\over
\d} }|+|\cos^{m'-n}{ \pi   {h\over
\d} } |\big)  \le  4(m-n)^{-\b'}. \cr} $$  

Further
$$ {2\over d\d}  \sum_{{1\le |j|<{d\over 2}\atop 1\le |h|<{\d\over 2}}\atop {\pi h\over \d} \in
I_{m-n}}   \Big|  \cos^{m-n}{ \pi   {h\over
\d} }-e^{ i\pi     {h\over \d}(m'-m)  }  \cos^{m'-n}{ \pi   {h\over
\d} }\Big| \le {4\over d\d} \sum_{{1\le |j|<{d\over 2}\atop 1\le |h|<{\d\over 2}}\atop {\pi h\over \d} \in
I_{m-n}}   1\le C\big({\log (m-n)\over m-n}\big)^{1/2}, $$ 
since  ${\pi h\over \d} \in
I_{m-n}$ means $h\le \d \pi \p_{m-n}$. We therefore get
 
$$\Big|{\bf \Delta} \big( (d,n), (\d, m)\big) -{\bf \Delta} \big( (d,n), (\d, m')\big)\Big|\le  C\big({\log (m-n)\over m-n}\big)^{1/2},$$ 
as claimed. 

  \bigskip\bigskip
It remains to compute $\rho_k(D)$. The lemma below is a classical tool. For the sake of completeness, we give a detailed proof.
 \medskip\noi{\gem  Lemma 2.15.} {\it Let $f\in \Z(X)$ and put $\rho_f(d)=\#\big\{0\le y<d: d|f(y)\big\}$. Then $\rho_f$ is a
multiplicative function.}
 
 \medskip\noi{\it Proof.}  Write $f(x)= a_0+a_1x+\ldots + a_nx^n$, $a_j\in \Z$, $0\le j\le n$. Let $d=d_1d_2$ with $(d_1,d_2)=1$. We first
establish
$\rho_f(d_1)\rho_f(d_2)\le
\rho_f(d )$. Let
$(y_1, y_2)$ be such that $0\le y_i<d_i$,   $d_i|f(y_i)$, $i=1,2$. There exists a unique integer $y$, $0\le y<d$ such that $y\equiv
y_i\, {\rm mod}(d_i)$,
$i=1,2$. Now, on writing $y=y_1+\ell d_1$
$$\eqalign{f(y)&= a_0+a_1(y_1+\ell d_1)+\ldots + a_n(y_1+\ell d_1)^n =a_0+(a_1 y_1+d_1 A_1  ) +\ldots + (a_n y_1^n +  d_1A_n)\cr & =
f(y_1)+ d_1 F ,\cr}$$
where $A_1, \ldots, A_n$, $F $ are integers. Thus $d_1|f(y)$. Similarly $d_2|f(y)$, and so $d|f(y)$. Now let $(y'_1, y'_2)$ be such that
$0\le y'_i<d_i$,  
$d_i|f(y'_i)$, $i=1,2$, and let $y'$ be the corresponding unique integer such that $0\le y'<d$, $y'\equiv
y'_i\, {\rm mod}(d_i)$, $i=1,2$, 
 and so $d|f(y')$. We have to prove the implication $y=y'\Rightarrow (y_1, y_2)=(y'_1, y'_2)$. But as $y=y'$,  we have
$$y= y_1+\ell d_1=y_2+k d_2=y'_1+\ell' d_1=y'_2+k' d_2.$$
And so $y_1-y'_1=  (\ell'-\ell) d_1$. Since $0\le y_1, y'_1<d_1$, this implies $y_1=y'_1$. Similarly $y_2=y'_2$, so that $(y_1,
y_2)=(y'_1, y'_2)$. Therefore $\rho_f(d_1)\rho_f(d_2)\le
\rho_f(d )$. 
\smallskip Conversely, let $0\le y<d$ be such that $d|f(y)$. Let $y_1, y_2$, $0\le y_i<d_i$ be such that  $y_i\equiv y \, {\rm
mod}(d_i)$, $i=1,2$. Then, in the same fashion
$$\eqalign{f(y_1)&= a_0+a_1(y +\ell d_1)+\ldots + a_n(y +\ell d_1)^n =a_0+(a_1 y +d_1 B_1  ) +\ldots + (a_n y ^n +  d_1B_n)\cr & =
f(y )+ d_1 G.\cr}$$
And so $d_1|f(y_1)$; similarly $d_2|f(y_2)$. Let   $0\le y'<d$ be such that $d|f(y')$, and let $(y'_1, y'_2)$ be the corresponding
pair of integers. Here again, we must prove the implication $  (y_1, y_2)=(y'_1, y'_2) \Rightarrow y=y'$. Write $y_1=y+\ell d_1$,
$y'_1=y'+\ell' d_1$. If $y_1=y'_1$, then $y-y'=(\ell'-\ell)d_1 $ so that $d_1|y-y'$. Similarly $y_2=y'_2$ implies $d_2|y-y'$. Thus
$d |y-y'$. As $0\le y,y'<d$ this implies that $y=y'$. Hence the implication $  (y_1, y_2)=(y'_1, y'_2) \Rightarrow y=y'$. And we deduce 
$ 
\rho_f(d )\le \rho_f(d_1)\rho_f(d_2) $. The proof is   complete. \cqfd

\medskip\noi{\gem  Proposition 2.16.} {\it We have}  
$$\rho_k(D)=\cases{\displaystyle{\ \ \ \ \ \ \, \prod_{p|D} p^{\lfloor {v_p(D)\over 2}\rfloor }} \qq &\hbox{ if $k=0$,}\cr
\displaystyle{ \prod_{ v_p(k)<{v_p(D)/
2}}(2p^{v_p(k)})\cdot\prod_{   v_p(k)\ge {v_p(D)/
2} }  p^{\lfloor{v_p(D)\over
2}\rfloor}} \qq &\hbox{ if $k\ge 1$.}\cr} $$
 In particular, if $D$ is squarefree, then $\rho_0(D)=1$.
 \medskip\noi{\it Proof.}   Let us first consider the   case: $1\le k\le m-n$, which is the
main case. In view of Lemma 2.15, it suffices to compute $\rho_k(p^r)$. Make a first observation: 
$$\rho_k(p^r)=2, \quad r=1,2,\ldots\qq\quad \hbox{if $p\not| k$}. $$ 

Indeed, if $(y,p)=1$, then $p^r|y+k$ and there is only one solution given by
$y\equiv -k \, {\rm mod}(p^r)$. If $y=p^s Y$, $(Y,p)=1$, $1\le s<r$, then $p^r| y(y+k)\Leftrightarrow p^{r-s} |(p^sY+k)$. And so $p|k$,
which was excluded. There is thus no solution of this kind. Finally, it remains one extra solution $y=p^r$. Thus $\rho_k(p^r)=2$.
 \smallskip
We thus concentrate on  the case $p|k$. We can range the solutions $y$ of the equation $p^r|y(y+k)$ in disjoint classes of type $y=p^sY$,
with
$(Y,p)=1$. When
$r=1,2$ or
$3$, there is a direct   computation and one find $$\rho_k(p^r)=\left\{\matrix { 
 1\ &{\rm if} \   r=1, 
 \cr
 p\ &{\rm if}\ r=2,   
  \cr
                                   2p  \ &{\rm if}\  r=3, \  v_p(k)=1. 
\cr
                                   p  \ &{\rm if}\  r=3, \  v_p(k)\ge 2. 
 }   
               \right.
 \eqno(2.40) $$
  Suppose now that $r\ge 4$ and put 
$$r' =\lfloor{r\over 2}\rfloor.$$  We have
$\rho_k(p^r)=\#\big\{1\le y\le p^r:
\  p^r| y(y+k) 
\big\}. 
$  
  If $(y,p)=1$, then $p^r| y(y+k)\Leftrightarrow p^r|y+k$ and so $p|y$, which is
excluded and there is no solution of this type. Apart from the trivial solution $y=p^r$, the other possible solutions are of type $y=p^s
Y$,  
$(Y,p)=1$,
$1\le s<r$; and we shall distinguish three cases:   
$$\hbox{{\it i}) $r'<s<r$, \qq{\it ii}) $s=r'$,\qq {\it iii}) $1\le s< r'$.}$$

\medskip\noi    {\it i}) Since $r'<s<r$, then $ r/2 \le  s $, and so $1\le r-s\le s$.  Further 
$p^r| y(y+k)$ means $p^{r-s}|  Y(p^sY+k)$ or
$p^{r-s}|  p^sY +k  $,   which is possible if and only if
$p^{r-s}|  k
$, namely $ r-s \le v_p(k)
$.  Thus 
$$\max(r'+1,r-v_p(k)) \le s<r.$$
 We have $Y\le p^{r-s}$, $(Y,p^{r-s})=1$. Their number
  is   $\phi(p^{r-s})$ where $\phi$ is Euler's function, and since $\phi(p^{r-s})=p^{r-s}(1-{1\over p}) $, 
the corresponding number of solutions is 
 $$\eqalign{  \sum_{\max(r'+1,r-v_p(k)) \le s\le r -1 } p^{r-s}(1-{1\over p})   &=   (1-{1\over p})\sum_{1\le  v  
\le (r-r'-1)\wedge v_p(k)}  
p^{v} \cr &  = p{p^{(\lfloor {r-1\over 2}\rfloor )\wedge v_p(k)}-1\over p-1}(1-{1\over p}) =  
p^{\lfloor {r-1\over 2}\rfloor\wedge v_p(k)} -1   . \cr}\eqno(2.41)$$

\medskip\noi  {\it ii})      We consider   solutions of type    $y=p^{r'}Y$, $(Y,p)=1$. 
\smallskip  ---  If $r$ is odd, $r=2r'+1$, then  $p^{2r'+1}| p^{ r' }Y(p^{ r' }Y+k)$ means
$p^{ r'+1}|   (p^{ r' }Y+k)$. So $p^{ r'}|k$ and thereby $v_p(k)\ge r'$.
  If $v_p(k)< r'$, there is thus no solution.  If $v_p(k)>  r'$, this implies that  $p|Y$ which is impossible and there is again no  
solution. 
 
The remainding case $v_p(k)=  r'$ will be the only one providing solutions. Write $k=p^{ r' }K $, $(K,p)=1$, then $p|Y+K$. Since
$(K,p)=1$, the solutions are the numbers
$Y$ such that
$1\le Y\le  p^{ r' +1}$ and 
$Y\equiv -K\, {\rm mod}(p)$. Let $1\le \kappa<p$ be such that $K  \equiv \kappa\, {\rm mod}(p)$.
The number of solutions is
$$\#\big\{Y\le  p^{ r' +1} : Y\equiv -K\, {\rm mod}(p)\big\} =\#\big\{ (p-\kappa) +jp: 0\le j<p^{ r' }\big\}= p^{ r'  }.\eqno(
2.42)$$

 \smallskip  --- If $r$ is even, $r=2r'$, then  $p^{2r' }| p^{ r' }Y(p^{ r' }Y+k)$ reduces to $p^{ r' }|  (p^{ r' }Y+k)$, so $p^{
r'}|k$. If $v_p(k)< r'$, there is no solution.  If $v_p(k)\ge   r'$,  write $k=p^{ r' }K$, $(K,p)=1$, this is always realized and   the
number of solutions is
$$\#\big\{Y\le  p^{ r' } :  (Y,p)=1\big\} =\phi(p^{ r' })=p^{ r' }(1-{1\over p}) .$$

 \medskip \noi {\it iii}) We consider the last type of  solutions: $y=p^sY$, $(Y,p)=1$,    $1\le s< 
r'$.  Notice first,  since $s<r'$ that  $s<r/2$, and so $r-s>r/2>s$.
As
$p^r| y(y+k)$ means 
 $p^{r-s}|   Yp^s+k $, we deduce that $p^s|k$, namely $ s\le v_p(k) $. 
 If $  v_p(k)<s $, there is no solution.  
  If $   v_p(k)>s $, then $p|Y$ which is impossible.  

If $s= v_p(k) $, which requires $v_p(k)<r'$, write $k=p^{v_p(k)}K$, $(K,p)=1$. Then
we get the equation $p^{r-2v_p(k)}|   Y +K $, so  $Y  
\equiv -K \, {\rm mod}(p^{r-2v_p(k)})$.  Notice that if $Y$ is a solution, then $(Y,p)= 1$, since $(K,p)=1 $. 
 Let $1\le \kappa<p^{r-2v_p(k)}$ be such that $K  \equiv \kappa\, {\rm mod}(p^{r-2v_p(k)})$.
The number of solutions is
 $$\eqalign{\#\big\{Y\le p^{r-v_p(k)}:    Y   \equiv -K \, {\rm mod}(p^{r-2v_p(k)})\big\}&  =\#\big\{ (p^{r-2v_p(k)}-\kappa)
+jp^{r-2v_p(k)}:\! 0\le j<p^{v_p(k)}\big\}\cr &=p^{v_p(k)}.\cr}\eqno(
2.43)$$
  
 Summarizing the case $r\ge 4$, if $r$ is odd, $r=2r'+1$ 
$$      \rho_k(p^r) =  \cases{    2p^{v_p(k)} & \quad\hbox{\rm if   $1\le v_p(k)\le r'$},
   \cr 
        p^{r'} & \quad\hbox{\rm if   $v_p(k)> r'$}.
   \cr
  }   
\eqno(2.44)   $$
And if $r$ is even, $r=2r'$  $$      \rho_k(p^r) =  \cases{     
   \cr
  2p^{v_p(k)} & \quad\hbox{\rm if  $1\le v_p(k)< r'$},
   \cr
 p^{r'} & \quad\hbox{\rm if   $v_p(k)\ge r'$}.
     }   
\eqno( 2.45)   $$
This remains true if $r=1,2,3$. Observe that ($v_p(k)< r'$, $r$ even) or ($v_p(k)\le  r'$, $r$ odd) are equivalent to $v_p(k)<{r\over
2}$.  Therefore, for
$r\ge 2$ 
 $$  \ \      \rho_k(p^r) =  \cases{    2  & \quad\hbox{\rm if   $ p\not| k$},
   \cr  2p^{  v_p(k)} & \quad\hbox{\rm if $v_p(k)<{r\over
2}$},
   \cr 
        p^{\lfloor{r\over
2}\rfloor}   \ & \quad\hbox{\rm if $v_p(k)\ge {r\over
2}$}.
     }   
\eqno(
2.46)   $$
 
 Consequently 
 $$\rho_k(D)= \prod_{p|D  }\rho_k(p^{  v_p(D) })=    \prod_{ v_p(k)<{v_p(D)/
2}}(2p^{v_p(k)})\cdot\prod_{   v_p(k)\ge {v_p(D)/
2} }  p^{\lfloor{v_p(D)\over
2}\rfloor} .
\eqno(
2.47)$$ \medskip

Finally,   consider the case  $k=0$, namely       $\rho_0(p^r)=\#\big\{1\le y\le p^r: \  p^r| y^2   \big\}  $. Notice that 
$\rho_0(p)=\#\big\{1\le y\le
p: \  p| y   \big\}=1 $. Let
$r>1$, and write   $y=p^sY$, $(Y,p)=1$ and $1\le s\le r$. If $s=r$, there is the trivial unique solution $y=p^r$. If $1\le s<r$,
 $p^r| y^2\Rightarrow 2s\ge r$, in which case, the number of solutions is  
$$\#\big\{Y\le p^{r- s}: (Y,p)=1\big\}=\phi(p^{r- s})=p^{r- s}(1-{1\over p}).\eqno(
2.48) $$
 Consequently
 $\rho_0(p^r)= 1+\sum_{r/2\le s< r}p^{r- s}(1-{1\over p})$ . 
If $r$ is even, $r=2r'$
$$\rho_0(p^r)= 1+\sum_{\s =1}^{r' } p^{\s}(1-{1\over p})= 1+p\sum_{u =0}^{r'-1 } p^{u}(1-{1\over p})= 1+p\Big({p^{r' }-1\over p-1}\Big)
 ( {p-1\over p})=p^{r' } ,$$
whereas if $r$ is odd, $r=2r'+1$,
 $\rho_0(p^r)= 1+\sum_{\s =1}^{r'} p^{\s}(1-{1\over p})=p^{r' }$.
Thus 
$$\rho_0(p^r)=p^{\lfloor {r\over 2}\rfloor }. $$
It follows that 
$$\rho_0(D)=\prod_{p|D} p^{\lfloor {v_p(D)\over 2}\rfloor }.\eqno(
2.49) $$
The proof is now complete. \cqfd 
We deduce
\medskip\noi{\gem Corollary 2.17.} {\it We have for any positive integer  $D$ and  any integer $k$,}  $$\rho_k(D)\le 2^{\o(D) }(k\wedge  
\sqrt D) ,\qq
\qq 
\rho_0(D)\le  
\sqrt D.
$$  
{\it Proof.} Immediate. 

\medskip\noi{\gem Corollary 2.18.} {\it We have for any positive integers $n,m$ and $D$,  
$$  \P \{ D|B_nB_m \}    \le   {1\over 2^{  m-n }\sqrt D
}+    {2^{\o(D)}( m-n)  \over D } + C_\e
 {     D  ^{{1/ 2} + \e}\over   \sqrt n }  .  
   $$
Further, for any $\e>0$, there exists a constant $C_\e$ depending on $\e$ only, such that} 
$$\P \{ 
d|B_n\, ,\, \d|B_m \}\le C   (m-n){ 2^{\o(d\d) } \over d\d
}+ C_\e  {     (d\d)  ^{(1 + \e)/2}\over  \sqrt n } .$$
{\it Proof.} By Corollary 2.17 $$\eqalign{   \sum_{k=0}^{m-n} { C^k_{m-n} \over 2^{ m-n }  } {\rho_k(D)\over D }&=  
{\rho_0(D)\over 2^{  m-n }D }+ 
\sum_{k=1}^{m-n} { C^k_{m-n}
\over 2^{  m-n }  } {\rho_k(D)\over D }  \le   {1\over 2^{  m-n }\sqrt D }+   \sum_{k=1}^{m-n} { C^k_{m-n}
\over 2^{  m-n }  } {2^{\o(D)} k\over D }\cr &  \le   {1\over 2^{  m-n }\sqrt D
}+    {2^{\o(D)}( m-n)  \over D } .\cr}$$
 
On using Proposition 2.19, we get  
 $$  \P \{ D|B_nB_m \}\le    \sum_{k=0}^{m-n} { C^k_{m-n} \over 2^{ (m-n)}  } {\rho_k(D)\over D }  + C_\e
\left({     D  ^{1 + \e}\over   n }\right)^{1/2 }     \le   {1\over 2^{  m-n }\sqrt D
}+    {2^{\o(D)}( m-n)  \over D } + C_\e
 {     D  ^{{1\over 2} + \e}\over   \sqrt n }  .  
   $$

  For proving the second estimate, notice first that $\P \{ 
d|B_n\, ,\, \d|B_m \}>0$ only if $m-n\ge (d,\d)$, since
 $\P \{ 
d|B_n\, ,\, \d|B_m \}  \le  \P \{ 
d|B_n\}\P \{ (d,\d)|B_{m -n} \} $.
    Now 
$$\eqalign{\P \{ 
d|B_n\, ,\, \d|B_m \} &  =\P \{ 
d|B_n\, ,\, \d|B_m\, ,\,B_m-B_n=0 \}+\P \{ 
d|B_n\, ,\, \d|B_m\, ,\,B_m-B_n>0 \}\cr &\le \P \{ 
d|B_n\, ,\, \d|B_n  \}+\P \big\{ 
d\d|\big(B_n^2 + B_n(B_m-B_n)\big) ,\,B_m-B_n>0 \big\}\cr &= \P \{ 
[d, \d]|B_n  \}+{1\over 2^n2^{m-n}}\sum_{k=1}^{m-n}C^k_{m-n}{1\over
d\d}\sum_{j=0}^{d\d-1}\sum_{h=0}^n   C^h_n e^{2i\pi{j\over d\d}  (h^2+kh)}.\cr}$$  
By (1.15), $\P \{ 
[d, \d]|B_n  \}\le C\big({1\over [d, \d] }+ {1\over \sqrt{n}}\big) \le
C\big({m-n \over  d  \d  }+ {1\over \sqrt{n}}\big)$. 
 And   by using (2.37) and Corollary 2.17
 $$\Big|{1\over 2^n2^{m-n}}\sum_{k=1}^{m-n}C^k_{m-n}{1\over
d\d}\sum_{j=0}^{d\d-1}\sum_{h=0}^n   C^h_n e^{2i\pi{j\over d\d}  (h^2+kh)} \Big|\le   C_\e {     (d\d)  ^{(1 + \e)/2}\over  \sqrt n }+{1\over  2^{m-n}}\sum_{k=1}^{m-n}{\rho_k(d\d)\over d\d
}$$ 
$$   \le  C_\e {     (d\d)  ^{(1 + \e)/2}\over  \sqrt n }+{1\over  2^{m-n}}\sum_{k=1}^{m-n}{ 2^{\o(d\d) }(k\wedge  
\sqrt{ d\d})\over d\d
} \le  C_\e {     (d\d)  ^{(1 + \e)/2}\over  \sqrt n }+ (m-n){ 2^{\o(d\d) } \over d\d
}.$$

Therefore
$$\eqalign{\P \{ 
d|B_n\, ,\, \d|B_m \}&\le C\big({m-n \over  d  \d  }+ {1\over \sqrt{n}}\big)+  (m-n){ 2^{\o(d\d) } \over d\d
}+C_\e  {     (d\d)  ^{(1 + \e)/2}\over  \sqrt n } \cr & \le C   (m-n){ 2^{\o(d\d) } \over d\d
}+ C_\e  {     (d\d)  ^{(1 + \e)/2}\over  \sqrt n }.\cr}$$
   \cqfd
 \noi {\it Proof of Proposition 2.10.} It suffices now to put together Propositions 2.14, 2.16 and Corollary 2.18.\cqfd
\medskip

\medskip We conclude this section by indicating another correlation estimate controlling the difference between ${\bf \Delta} \big(
(d,n), (\d, m)\big)$ and ${\bf \Delta} \big( (d,n), (\d, m')$  once $m,m'$ are not too close to $n$.

\medskip 
\noi{\gem Proposition 2.19.} {\it For some $n_0 $   sufficiently large, and 
  $m'\ge m \ge n+n_0$, }
$$\Big|{\bf \Delta} \big( (d,n), (\d, m)\big) -{\bf \Delta} \big( (d,n), (\d, m')\big)\Big|\le  C\big({\log (m-n)\over m-n}\big)^{1/2}.$$ 
 
\medskip\noi {\it Proof.} In view of (2.12)
$$  {\bf \Delta} \big( (d,n) , (\d, m)\big)  
 ={1\over d\d}\sum_{1\le |h|<\d/2\atop1\le |j|<d/2}    e^{ i\pi  (  {j\over d}n +   {h\over \d}m) }\cos^{m-n}{ \pi  
{h\over
\d} }  \Big\{  \cos^n{ \pi  (  {j\over d} +  
{h\over \d})   }  -    \cos^{n} {
\pi j\over d}  \cos^{n} {
\pi h\over \d}\Big\}.   $$  
So that for $n\le m\le m'$,
$$\displaylines{{\bf \Delta} \big( (d,n), (\d, m)\big) -{\bf \Delta} \big( (d,n), (\d, m')\big)=  {1\over d\d}\!\sum_{1\le
|h|<\d/2\atop1\le |j|<d/2}    \Big\{\cos^n{ \pi  (  {j\over d} +   {h\over \d})   }-\cos^{n} {
\pi j\over d}     \cos^{n} {
\pi h\over \d}\Big\}\hfill\cr\hfill\times \ e^{ i\pi     {j\over d}n   }\bigg[e^{ i\pi     {h\over \d}m  }  \cos^{m-n}{ \pi  
{h\over
\d} }-e^{ i\pi     {h\over \d}m'  }  \cos^{m'-n}{ \pi   {h\over
\d} }\bigg].\cr}$$
 Then  we may bound this difference as follows:
$$\Big|{\bf \Delta} \big( (d,n), (\d, m)\big) -{\bf \Delta} \big( (d,n), (\d, m')\big)\Big|\le   {2\over d\d}\!\!\sum_{1\le
|h|<\d/2\atop1\le |j|<d/2}   \Big|  \cos^{m-n}{ \pi   {h\over
\d} }-e^{ i\pi     {h\over \d}(m'-m)  }  \cos^{m'-n}{ \pi   {h\over
\d} }\Big|.$$ 
 
In view of (2.13) and Case I of the proof of Proposition 2.1, if ${\pi h\over \d} \in I'_{m-n}$, then 
    $|\cos {\pi h\over \d}|\le \cos \p_{m-n}$. And so, for some $n_0 $   sufficiently large, and 
  $m-n\ge n_0$, 
  $$|\cos {\pi h\over \d}|^{m'-n} \le|\cos {\pi h\over \d}|^{m-n} \le  \cos^{m-n} \p_{m-n}    
=e^{-2(m-n)\sin^2(\p_{m-n}/2)}\le (m-n)^{-\b'}  .$$ 
    
 Hence,  
$$\eqalign{ {2\over d\d}  \sum_{{1\le |j|<{d\over 2}\atop 1\le |h|<{\d\over 2}}\atop {\pi h\over \d} \in
I'_{m-n}}   \Big|  \cos^{m-n}{ \pi   {h\over
\d} }-&e^{ i\pi     {h\over \d}(m'-m)  }   \cos^{m'-n}{ \pi   {h\over
\d} }\Big|\cr & \le {2\over d\d} \sum_{{1\le |j|<{d\over 2}\atop 1\le |h|<{\d\over 2}}\atop {\pi h\over \d} \in
I'_{m-n}}   \big(  |  \cos^{m-n}{ \pi   {h\over
\d} }|+|\cos^{m'-n}{ \pi   {h\over
\d} } |\big)  \le  4(m-n)^{-\b'}. \cr} $$  

Further
$$ {2\over d\d}  \sum_{{1\le |j|<{d\over 2}\atop 1\le |h|<{\d\over 2}}\atop {\pi h\over \d} \in
I_{m-n}}   \Big|  \cos^{m-n}{ \pi   {h\over
\d} }-e^{ i\pi     {h\over \d}(m'-m)  }  \cos^{m'-n}{ \pi   {h\over
\d} }\Big| \le {4\over d\d} \sum_{{1\le |j|<{d\over 2}\atop 1\le |h|<{\d\over 2}}\atop {\pi h\over \d} \in
I_{m-n}}   1\le C\big({\log (m-n)\over m-n}\big)^{1/2}, $$ 
since  ${\pi h\over \d} \in
I_{m-n}$ means $h\le \d \pi \p_{m-n}$. We therefore get
 
$$\Big|{\bf \Delta} \big( (d,n), (\d, m)\big) -{\bf \Delta} \big( (d,n), (\d, m')\big)\Big|\le  C\big({\log (m-n)\over m-n}\big)^{1/2},$$ 
as claimed.\cqfd


\medskip  \medskip\noi\centerline{\gum 3. Increments of sums of divisors of Bernoulli sums.}
\medskip \medskip\noi  
Let $0<\t<1/6$.   Put for any positive integer $n$
 $$H_n=H_n({\cal D}, \t) = \sum_{ d\le n^\t, \, d\in {\cal D} }\big( {\bf 1}_{d|B_n}- \P\{d|B_n\}\big).  $$
Let also  an increasing sequence ${\cal N} $  satisfying the growth condition ${\cal G}_\rho$  for some 
$\rho>0$. Let $\eta>0$
  and put 
$$\widetilde H_n=\widetilde H_n ({\cal D} ) = \sum_{ d<\eta \sqrt{n\over \log n}, \,  d\in {\cal D} }\big( {\bf 1}_{d|B_n}-
\P\{d|B_n\}\big).  $$
    In this section  we   establish the following result. 
 \medskip \noi{\gem Theorem 3.1.} {\it a)  For any $\e>0$, there exist    constants $C_\e$ and $i_\e$, such that for every
$i$ and $j$ with $\min(i,j-i)\ge i_\e$, 
  $$ \E\big(\sum_{i\le n\le j} H_n\big)^2\le  C_\e\sum_{i\le n\le j}n^{\e} .   $$
b) There exist   constants $\eta_0>0$, $C<\infty$ such that for $\eta\le \eta_0$ and $j\ge i$}
 $$ \E\big(\sum_{i\le n\le j\atop n\in{\cal N}} \widetilde H_n\big)^2\le  C\sum_{i\le n\le j\atop n\in {\cal N}}(\log n)^4 .   $$

\noi{\it Proof.}    We   rewrite   $ \E\big(\sum_{i\le n\le j} H_n\big)^2$ as follows
$$ \E\big(\sum_{i\le n\le j} H_n\big)^2=  \sum_{i\le n\le j} \E H_n^2+  2\sum_{i\le n\le j}\sum_{  n<m\le j} \E H_nH_m:= A+2B.\eqno(3.1)$$
For the other increment $ \E\big(\sum_{i\le n\le j} \widetilde H_n\big)^2$ we operate identically.  Let $   0<c<1/5$   and choose  $H=4c$. 
We split the sum $B$   into two subsums as follows:
$$\eqalign{ B&= \sum_{i\le n\le j}\sum_{  n<m\le j}\sum_{ d\le n^\t,\d\le m^\t\atop d,\d\in {\cal D} } {\bf \Delta} \big( (d,n), (\d,
m)\big) =\sum_{i\le n\le j}\sum_{  n<m\le n+n^{H}  } \sum_{ d\le n^\t,\d\le m^\t\atop d,\d\in {\cal D} }{\bf \Delta} \big( (d,n),
(\d, m)\big)\cr &\quad+\sum_{i\le n\le j}\sum_{  n+n^{H} <m\le j} \sum_{ d\le n^\t,\d\le m^\t\atop d,\d\in {\cal D} }{\bf \Delta} \big(
(d,n), (\d, m)\big)
 :=B_1+B_2.\cr}\eqno(3.2)$$
 
The sum $B_1$ is really typical from the "small increments" case. And we will see that this sum, which will be examined by means of
Proposition 2.10, produces the strongest contribution. Concerning the   sum $B_2$, let $\d_1>2$ arbitrary but fixed. By Propositions 2.5
and 2.9, we know that  there exist  constants 
$n_0$, $C $ such that for any $  n\ge n_0$, 
 $$\line{$  \displaystyle{\sup_{d<{\pi\over \d_1\sqrt{ \a}}  \sqrt{n\over \log n}  \atop \d<{\pi\over \d_1\sqrt{ \a}} \sqrt{m\over\log m} 
}}\big|{\bf \Delta}
\big( (d,n), (\d, m)\big)\big|
 \le 
    \left\{\matrix{    C  n  ({  \log (m-n)\over m-n})^{ 1/c} 
     \qq\ \   &{\rm if}\   n+n^c \le m\le 2n ,
\cr &\cr    C  
n({  \log n\over n})^{{1/ 2c}}({ \log (m-n)\over
m-n})^{1/2} &      {\rm if}\    m\ge 2n .\qq\qq\cr
 }  \right.  \hfill$}\eqno(3.3)$$

  \medskip\noi     {\bf (1) Estimating the sum $\bf B_2$.} We claim that
$$\sum_{i\le n\le j}\sum_{  n+n^{H} <m\le j} \sum_{d<{\pi\over \d_1\sqrt{ \a}}  \sqrt{n\over \log n}\atop\d<{\pi\over \d_1\sqrt{ \a}} \sqrt{m\over\log m}
}{\bf \Delta} \big( (d,n), (\d, m)\big)\le  C(j-i).\eqno(3.4)$$
 \noi --- If $j\le 2i$, using (3.3) we get
   $$\eqalign{\sum_{i\le n\le j}\sum_{  n+n^{H}<m\le
j  }\sum_{d<{\pi\over \d_1\sqrt{ \a}}  \sqrt{n\over \log n} \atop\d<{\pi\over \d_1\sqrt{ \a}} \sqrt{m\over\log m} } n ({  \log (m-n)\over
m-n})^{ 1/c} 
   &\le C\sum_{i\le n\le j}n^{2 - H/ c }(\log  n)^{1/c-2}\sum_{  n+n^{H}<m\le
j  }   1     \cr
&\le C(j-i)\sum_{i\le n\le j}n^{-2}\log^{ 1/c-2} n 
  \le C(j-i). \cr} $$
  --- And if $j\ge 2i$, 
    $$\eqalign{\sum_{i\le n\le j}\sum_{ \min(2n,j) <m\le
j  }&\sum_{d<{\pi\over \d_1\sqrt{ \a}}  \sqrt{n\over \log n} \atop\d<{\pi\over \d_1\sqrt{ \a}} \sqrt{m\over\log m} } n     ({\log  n\over n})^{{1/
2c}}  ({
\log (m-n)\over m-n})^{1/2} \cr &\le  C\sum_{i\le n\le j}n^{3/2-{1/
2c}}    ({\log  n })^{{1/
2c}-1/2} \sum_{ \min(2n,j) <m\le
j  }  ({m
 \over \log m})^{1/2} ({
\log m\over m})^{1/2}     \cr &\le  C \sum_{  n\ge 1}n^{3/2-{1/
2c}}    ({\log  n })^{{1/
2c}-1/2} \sum_{ \min(2n,j) <m\le
j  } 1  
 \le C (j-i), \cr} $$
  since $c<1/5$.
 
\medskip  {\bf Remark.} --- Let $0<\t<1/2$. It also follows from (3.4) that 
$$
\sum_{i\le n\le j}\sum_{  n+n^{H} <m\le j} \sum_{ d\le n^\t,\d\le m^\t\atop d,\d\in {\cal D} }{\bf \Delta} \big( (d,n), (\d,
m)\big) \le C_\t (j-i).\eqno(3.5a)$$
--- Let ${\cal N}$ be an increasing sequence of integers. It also follows trivially from (3.4) that
$$\sum_{i\le n\le j\atop n\in {\cal N}}\sum_{  n+n^{H} <m\le j\atop m\in {\cal N}} 
\sum_{{d<{\pi\over \d_1\sqrt{ \a}}  \sqrt{n\over \log n}\atop\d<{\pi\over
\d_1\sqrt{ \a}}
\sqrt{m\over\log m}}\atop d,\d\in {\cal D} }{\bf \Delta} \big( (d,n), (\d, m)\big)\le   C\sum_{i\le n\le j\atop n \in {\cal N}}1  
.\eqno(3.5b)$$

\bigskip\noi     {\bf (2) Estimating the sum $\bf B_1$.}  Let $h>0$ be some small number. 
 By Proposition 2.10,  
  $$ \P \{ 
d|B_n\, ,\, \d|B_m \}   \le C    { (m-n)2^{\o(d\d) } \over d\d
}+C_h  {     (d\d)  ^{(1 + h)/2}\over  \sqrt n } . $$
Thus  
$$\eqalign{ B_1  &=\sum_{i\le n\le j\atop  n<m\le n+n^{H}  } \sum_{ d\le n^\t,\d\le m^\t\atop d,\d\in {\cal D} }{\bf \Delta} \big(
(d,n), (\d, m)\big)   \le \sum_{i\le n\le j\atop   n<m\le n+n^{H}  } \sum_{ d\le n^\t,\d\le m^\t\atop d,\d\in {\cal D} }\P \big\{ d|
B_n\, ,\, \d|B_m \big\}
\cr & \le \sum_{i\le n\le j\atop  n<m\le n+n^{H}  } \sum_{ d\le n^\t,\d\le m^\t\atop d,\d\in {\cal D} }
\Big( C   (m-n){ 2^{\o(d\d) } \over d\d
} + C_h
 {     (d\d)  ^{{1\over 2} + h}\over   \sqrt n }\Big) := B_1^1+B_1^2 . \cr}\eqno(3.6)$$
For the sum $B_1^2$  we have since $\t<1/6$
$$B_1^2=\sum_{i\le n\le j}\sum_{  n<m\le n+n^{H}  } \sum_{ d\le n^\t,\d\le m^\t\atop d,\d\in {\cal D} } 
  {     (d\d)  ^{(1 + \e)/2}\over   n^{1/2 } }  \le  \sum_{i\le n\le j}n^{ \t(3+\e)+H-1/2}\le C_\t(j-i),\eqno(3.7)$$
if $\e$, $c$, ($H=4c$) are small enough, which we do assume. 
 For the sum $B_1^1$, we have  
 $$ \sum_{ d\le n^\t,\d\le
m^\t\atop d,\d\in {\cal D} }    {2^{\o(d\d)} \over d\d }       \le 2 \sum_{ D\le (nm)^\t  } \#\big\{  d\le n^\t,\d\le
m^\t: D=d\d\big\}\cdot
  {2^{\o(D)}  \over D } 
 \eqno(3.8) $$ 
But     $D=d\d$ occurs, given $d$, only for one choice of $\d$: $\d=D/d$. Further, the number of possible $d$ cannot exceed  the number of divisors of
$D$: $d(D)$. Thus $\#\big\{  d\le n^\t,\d\le
m^\t: D=d\d\big\}\le  \#\big\{  d ,\d : D=d\d\big\}\le d(D)$. And it is well-known that $d(N)={\cal O}_\e(N^\e)$. Thereby $d(D)={\cal O}_\e(D^\e) $,  for
$D$ large, say $D\ge \D_\e$. And so, if $D\ge
\D_c$  
$$   d(D)\le C_c D^c\le C_c n^c.$$
Obviously 
$$   \sum_{    D<\D_c } \#\big\{  d\le n^\t,\d\le
m^\t: D=d\d\big\}\cdot
  {2^{\o(D)}  \over D } \le \sum_{    D<\D_c }d(D)\cdot
  {2^{\o(D)}  \over D } \le K(\D_c)  
 , $$
whereas  $$   \sum_{\D_c\le  D\le (nm)^\t  } \#\big\{  d\le n^\t,\d\le
m^\t: D=d\d\big\}\cdot
  {2^{\o(D)}  \over D } \le   C_c n^c\sum_{ D\le 2 n ^{2\t}  } 
  {2^{\o(D)}  \over D }
 , $$
 
  Put for a while 
 $ F(x)=\sum_{ k\le x }  {2^{\o(k)}\over  k} $, and recall ([T] p. 60 Exercise 5) that 
$$F(x) ={C_0\over 8\log 2} x(\log x)^2 +{\cal O}(x\log x), \qq C_0=\prod_{p>2}(1+{1\over p(p-2) }) .$$ 
Using Abel summation, we   deduce  that 
$$\sum_{  2\le D\le N}{2^{\o(D)}\over  D }=-{F(1)\over 2}+\sum_{j=2}^{N-1} {F(j)\over j(j+1)}+ {F(N)\over N}\le C\sum_{j=1}^N{(\log
j)^2\over j}\le C(\log N)^3. \eqno(3.9)$$
 
Hence  
  $$ \sum_{ d\le n^\t,\d\le m^\t\atop d,\d\in {\cal D} }  {2^{\o(d\d)}\over  d\d } \le \sum_{ D\le (nm)^\t  }   
{2^{\o(D)}\over  D} \le  C_\t(\log nm)^3.\eqno(3.10) $$
Summarizing, for $n$ large
$$\eqalign{ B_1^1&=\sum_{i\le n\le j}\sum_{  n<m\le n+n^{H}  } \sum_{ d\le n^\t,\d\le m^\t\atop d,\d\in {\cal D} }  {2^{\o(d\d)}( m-n) 
\over d\d } \le  C_{\t,   c}  
\sum_{i\le n\le j}n^c\sum_{  n<m\le n+n^{H}  }(m-n) (\log nm)^3\cr &\le  
C_{\t,   c}   \sum_{i\le n\le j}n^c(\log 2n^2 )^3\sum_{  n<m\le n+n^{H}  }(m-n)\le  
C_{\t,   c}  \sum_{i\le n\le j} n^{2H+c}(\log 2n^2 )^3 . \cr}\eqno(3.11) $$ 
 Thereby   for $n$ large 
$$\eqalign{B_1  
  \le C_{\t,   c} \sum_{i\le n\le j}  n^{12c} .\cr}\eqno(3.12)$$

  \noi    {\bf (3) Estimating the sum $\bf A$.}  Now we turn to the sum
  $A=\sum_{n=i}^j 
  \E H_n^2 , $
and begin with 
$$\E H_n^2 = \sum_{    
d,\d\in {\cal D} \atop  d,\d\le      n ^\t   } \P\{[d,\d]|B_{n}\}-\P\{d|B_{n}\}\P\{\d|B_{n}\}.  $$
In view of (1.16)
$$\sup_{    
u\in {\cal D} \atop  u\le      n ^\t   }\big|\P\big\{u|B_n\big\}- {1\over u}  \big|\le C ({u \over n})\le
Cn^{-(1-\t)}. $$
 so that 
$$\Big| \P\{[d,\d]|B_{n}\}-\P\{d|B_{n}\}\P\{\d|B_{n}\}- \big({1\over [d,\d]}-{1\over  d \d }\big)\Big|\le C
n^{-(1-2\t)}.\eqno(3.13)$$
This estimate is efficient only if ${1\over [d,\d]}-{1\over  d \d }$ is small. If $d,\d$ are coprimes the latter quantity vanishes and
(3.13) makes sense. Otherwise the correct order of ${1\over [d,\d]}-{1\over  d \d }$  is given by ${1\over [d,\d]} $, and is for $\t$
small, much bigger than
$n^{-(1-2\t)}$. And then, one has advantage to use the simple bound (see (1.19))
$$\big| \P\{[d,\d]|B_{n}\}-\P\{d|B_{n}\}\P\{\d|B_{n}\} \big|\le   {2\over [d,\d]}.\eqno(3.14)$$
According to  Eq. 18.2.1 p.263 of [HW]     
and   Eq. (B) p.81 of [R]  (see   [Wi] for a proof)  we recall  that  $\sum_{n=1}^N d^2(n) 
\sim  ({N\over \pi^2}) \log^3 N$. Thus 
$$0\le \E H_n^2  \le  2\sum_{    
d,\d\in {\cal D} \atop  d,\d\le      n^\t   } {1\over [d,\d]}\le C\sum_{    
h \le       n^{2\t}   } { d^2(h)\over h}\le C\log^4 n,\eqno(3.15)$$
where we used Abel summation for obtaining the last inequality. Thereby,
$$A= \sum_{n=i}^j 
  \E H_n^2 \le C\sum_{n=i}^j   \log^4 n. \eqno(3.16)$$
Combining (3.5), (3.11) with (3.18)  shows that there exists a constant $C$ depending on $c$, such that for $i$ large enough, say $i\ge i_c$,  
$$\E\big(\sum_{i\le n\le j} H_n\big)^2\le C_\t\sum_{i\le n\le j}(\log  n  )^3 n^{2H}\le C_{\t, c}\sum_{i\le n\le j}  n^{13c} .
\eqno(3.17)$$

Now we estimate the other increment. The major difference in comparison with the
above   lies in the fact that the sum $B_1$ disappears, once $4c<\rho$ which we do assume. We start similarly to (3.1) with 
 $$ \E\Big(\sum_{i\le n\le j\atop n\in {\cal N}} H_n\Big)^2=  \sum_{i\le n\le j\atop n\in {\cal N}} \E H_n^2+  2\sum_{i\le n\le
j\atop n\in {\cal N}}\sum_{  n<m\le j\atop m\in {\cal N}}
\E H_nH_m:= A'+2B', $$
 where 
$$\eqalign{ B'&= \sum_{i\le n\le j\atop n\in {\cal N}}\sum_{  n<m\le j\atop m\in {\cal N}}
\sum_{ {d<\eta  \sqrt{n\over \log n} \atop\d< \eta
\sqrt{m\over\log m}}\atop d,\d\in {\cal D} } {\bf \Delta}
\big( (d,n), (\d, m)\big) =\sum_{i\le n\le j\atop n\in {\cal N}}\sum_{  n+n^{H} <m\le j\atop n\in {\cal N}} \sum_{ {d<\eta 
\sqrt{n\over
\log n} \atop\d< \eta
\sqrt{m\over\log m}}\atop d,\d\in {\cal D} }{\bf
\Delta} \big( (d,n), (\d, m)\big)
.\cr} $$ 
 If $\eta\le \d_1\sqrt{ \a}$, by (3.5b) we get   
$$ B'\le C  \sum_{i\le n\le j\atop n\in {\cal N}}1.\eqno(3.18) $$
  And in a same fashion as for getting (3.16) 
$$A'\le C \sum_{i\le n\le j\atop n \in {\cal N}} 
     \log^4 n. \eqno(3.19)$$
Finally 
$$ \E\big(\sum_{i\le n\le j\atop n \in {\cal N}} H_n\big)^2\le  C\sum_{i\le n\le j\atop n \in {\cal N}} \log^4 n .  \eqno(3.20)$$
 The proof is now complete.  \cqfd


\medskip\noi\centerline{\gem 4. Growth of sums of divisors of Bernoulli sums.}
 \medskip\noi
In this section, we  prove the main results of the paper. We begin with recalling a useful convergence result of
G\'al-Koksma type.

 \medskip\noi{\gem Lemma 4.1.} ([We5], Corollary 1.5) {\it  Let the random variables $\xi= \{\xi_i, i\ge 1\}$ satisfy
the following
assumption:
$$\E\big| \sum_{\ell=i}^j \xi_\ell\big|^2\le  \sum_{\ell=i}^j m_\ell , \qq \qq (i\le j) \eqno  (4.1)   $$
 where $  \{m_\ell, i\ge 1\}$ is a sequence of non negative reals such that the
series $ \sum_{\ell=1}^\infty u_\ell$ diverges. Put $M_n=\sum_{\ell=1}^nm_\ell$. Assume that  $$\log{M_n\over m_n}={\cal O}\big(\log
M_n\big).\eqno(4.2)$$

\noindent Then for any $\tau>1$,} 
$$\sum_{1\le \ell\le n}
\xi_\ell\buildrel{a.s.}\over={\cal O}_\tau\Big(M_n^{1/2} \bigl(\log (1+
M_n)\bigr)^{1+\tau/2}\Big)   
\eqno(4.3)$$
  \medskip\noi{\it Proof of Theorem 1.1.} By Theorem 3.1, for any $\e>0$ and  
$i$, $j$ such that  $\min(i,j-i) $ is large enough  
  $$ \E\big(\sum_{i\le n\le j} H_n\big)^2\le  C_\e\sum_{i\le n\le j}n^{\e} .   $$
 Thus condition (4.1) is fulfilled with $m_\ell= \ell^\e$. Further condition (4.2) trivially holds.   We also notice that 
$$\eqalign{M_{   {\cal D}}(n) &=\sum_{k=1}^n \E d_{\cal D} (B_{k}   ) =\sum_{k=1}^n\sum_{d\in {\cal D}\atop d\le k^\t}   \P\{ d |B_{k} 
\} \ge C\sum_{k=1}^n\sum_{d\in {\cal D}\atop d\le k^\t}   {1\over d}\ge C\sum_{d\in {\cal D}\atop d\le k^\t} \sum_{k=\lfloor
d^{1/\t}\rfloor}^n   {1\over d} 
\cr & =
C   \sum_{d\in {\cal D}\atop d\le n^\t} { (n - \lfloor d^{1/\t}\rfloor)\over d}\ge C    \sum_{d\in {\cal D}\atop d\le (n /2)^\t}  {(n -
\lfloor d^{1/\t}\rfloor)\over d} \ge Cn  \sum_{d\in {\cal D}\atop d\le (n /2)^\t}\gg n.
\cr}$$ Thus Theorem 1.1 is now  a direct consequence of Lemma 4.1.\cqfd
 
 \medskip\noi{\it Proof of Theorem 1.2.}  By theorem 3.1, there exist   constants $\eta_0>0$, $C<\infty$ such that for $\eta\le \eta_0$
and
$j\ge i$ 
 $$ \E\big(\sum_{i\le \nu_k\le j } \widetilde H_{\nu_k}\big)^2\le  C\sum_{i\le \nu_k\le j }(\log \nu_k)^4 .   $$
Put 
  $M_J= \sum_{    \ell\le J} 
    m_\ell $,  $m_\ell= \log^4 \nu_\ell$. 
 Condition (4.2) of  Lemma 4.1 is further trivially satisfied.  
 Then for any $b>3/2$,  
$$  \sum_{\ell=1}^J d_{\eta,\cal D} (B_{\nu_\ell}   )  \buildrel{a.s.}\over{=}M_{\eta,{\cal N},  {\cal D}}(J) +{\cal O}_\e\Big( 
M_J^{1/2} 
\log^bM_J\Big)    .
\eqno(4.4)$$ 
 Now  by (1.12), $N$ large enough so that there is $n\in {\cal N}$ such that $\eta \sqrt{n\over \log n}\ge \min\{{\cal D}\}$
 $$\eqalign{M_{\eta,{\cal N},  {\cal D}}(N) &  \ge
C\sum_{  n\le N   \atop n\in{\cal N}}\sum_{ d<\eta \sqrt{n\over \log n}\atop   d\in {\cal D} } {1\over d} \ge C\sum_{  n\le N   \atop
n\in{\cal N}}1
 \cr}\eqno(4.5)$$ 
Hence  if ${\cal N}$ grows at most polynomially, letting $N=\nu_J$ we get
$$M_J\le J\log^4 \nu_J\le CJ\log^4  J\le  M_{{\cal N},  {\cal D}}(J)\log^4 M_{{\cal N},  {\cal D}}(J).\eqno(4.6)$$
And in this case,
for any $b>7/2$,  
$$  \sum_{\ell=1}^J d_{\cal D} (B_{\nu_\ell}   )  \buildrel{a.s.}\over{=}M_{{\cal N},  {\cal D}}(J) +{\cal O}_\e\Big( 
 M_{{\cal N},  {\cal D}}(J)^{1/2} 
\log^b M_{{\cal N},  {\cal D}}(J)\Big)    .
\eqno(4.7)$$


\noi {\gum References}
 \medskip   \vskip 1pt\noi [BW]  {Berkes I., Weber M.}   [2007]  {\sl A law of the
iterated logarithm for arithmetic functions},   Proc. of the Amer. Math. Soc. {\bf 135} no4, 1223-1232.
 \vskip 1pt \noi [HR] {Halberstam H., Richert H. E.} [1974]   {\sl Sieve Methods}, Academic Press.
  \vskip 1pt \noindent[H] {Hardy G.H.},
[1963]  {\sl Divergent series}, Oxford at the Clarendon Press.
  \vskip 1pt  \noi [HW] {Hardy G.H., Wright E. M.}, [1979]  {\sl An
 introduction to the theory of numbers}, Oxford at the Clarendon
 Press, Fifth ed.
  \vskip 1pt  \noi[S1] {S\'ark\H ozy A.}, [1978] {\sl On difference sets of sequence of integers}, Acta Math. Acad. Sci. Hungar. {\bf 31}
(1-2), 125-149.
 \vskip 1pt  \noi[S2] {S\'ark\H ozy A.}, private communication.   
 \vskip 1pt  \noi [T] {Tenenbaum G.}, [1990]: {\sl Introduction  \`a
la  th\'eorie analytique et probabiliste des nombres},
Revue de l'Institut Elie
Cartan  {\bf 13}, D\'epartement de Math\'ematiques de
l'Universit\'e de Nancy I. \vskip 1pt\noi [We1]  {Weber M.}   [2004]   {\sl An arithmetical property of
Rademacher sums}.  Indagationes Math.  {\bf 15}, No.1, 133-150.
 \vskip 1pt\noi [We2]  {Weber M.}   [2005]  {\sl Divisors, spin sums  and
the functional equation of the Zeta-Riemann function},   
  Periodica Math. Hungar.  {\bf 51} (1),  1-13.
  \vskip 1pt \noindent [We3]  {Weber M.}    [2005]\ {\sl Small divisors of Bernoulli sums},  Indag. Math. {\bf 18  (2)}, (2007), p.281--293.
   \vskip 1pt\noi [We4]   {Weber M.},  [2006] {\sl On the order of magnitude
of the divisor function},   Acta Math. Sinica {\bf 22}, No2, 377-382.
\vskip 1pt\noi [We5]    Weber M.  [2006] {\sl Uniform bounds under increment conditions}.    Trans. Amer. Math. Soc. {\bf 358} 
no.2, 911-936. 
  \vskip 1pt \noindent [Wi]  {Wilson B. M.}  [1922]. {\sl Proofs of some formulae enunciated by Ramanujan}, Proc. of the London Math. Soc.
(2) {\bf 21},    235-255.
   \medskip 
 \noi {\ph Address} : U.F.R. de Math\'ematique (IRMA), Universit\'e Louis-Pasteur et C.N.R.S., 7 rue Ren\'e
Descartes,  F-67084 Strasbourg Cedex  
 
\noi  {\ph Email} :  \tt  weber@math.u-strasbg.fr \par

\bye